# Gradient descent with nonconvex constraints: local concavity determines convergence

Rina Foygel Barber and Wooseok Ha

10.17.17

### Abstract

Many problems in high-dimensional statistics and optimization involve minimization over nonconvex constraints—for instance, a rank constraint for a matrix estimation problem—but little is known about the theoretical properties of such optimization problems for a general nonconvex constraint set. In this paper we study the interplay between the geometric properties of the constraint set and the convergence behavior of gradient descent for minimization over this set. We develop the notion of local concavity coefficients of the constraint set, measuring the extent to which convexity is violated, which govern the behavior of projected gradient descent over this set. We demonstrate the versatility of these concavity coefficients by computing them for a range of problems in low-rank estimation, sparse estimation, and other examples. Through our understanding of the role of these geometric properties in optimization, we then provide a convergence analysis when projections are calculated only approximately, leading to a more efficient method for projected gradient descent in low-rank estimation problems. Projected gradient descent; local concavity; prox-regular sets; low-rank estimation

## 1 Introduction

Nonconvex optimization problems arise naturally in many areas of high-dimensional statistics and data analysis, and pose particular difficulty due to the possibility of becoming trapped in a local minimum or failing to converge. Nonetheless, recent results have begun to extend some of the broad convergence guarantees that have been achieved in the literature on convex optimization, into a nonconvex setting.

In this work, we consider a general question: when minimizing a function $\mathbf{g}(x)$ over a nonconvex constraint set $\mathcal{C} \subset \mathbb{R}^d$,

$$\widehat{x} = \arg\min_{x \in \mathcal{C}} \mathbf{g}(x),$$

what types of conditions on $\mathbf{g}$ and on $\mathcal{C}$ are sufficient to guarantee the success of projected gradient descent? More concretely, when can we expect that optimization of this nonconvex problem will converge at essentially the same rate as a convex problem?

In examining this question, we find that local geometric properties of the nonconvex constraint set $\mathcal{C}$ are closely tied to the behavior of gradient descent methods, and the main results of this paper study the equivalence between local geometric conditions on the boundary of $\mathcal{C}$, and the local behavior of optimization problems constrained to $\mathcal{C}$.

The main contributions of this paper are:

- We develop the notion of *local concavity coefficients* of a nonconvex constraint set $\mathcal{C}$, characterizing the extent to which $\mathcal{C}$ is nonconvex relative to each of its points. These coefficients, a generalization of the notions of *prox-regular sets* and *sets of positive reach* in the analysis literature, bound the set's violations of four different characterizations of convexity—e.g. convex combinations of points must lie in the set, and the first-order optimality conditions for minimization over the set—with respect to a structured norm, such as the $\ell_1$ norm for sparse problems, chosen to capture the natural structure of the problem. The local concavity coefficients allow us to characterize the geometric properties of the constraint set $\mathcal{C}$ that are favorable for analyzing the convergence of projected gradient descent. Our key results Theorems 1 and 2 prove that these multiple notions of nonconvexity are in fact exactly equivalent, shedding light on the interplay between geometric properties such as curvature, and optimality properties such as the first-order conditions, in a nonconvex setting.

- We next prove convergence results for projected gradient descent over a nonconvex constraint set, minimizing a function g assumed to exhibit restricted strong convexity and restricted smoothness. We also allow for the projection step, i.e. projection to $\mathcal{C}$, to be calculated approximately, which enables greater computational efficiency. Our main convergence analysis shows that, as long as we initialize at a point $x_0$ that is not too far away from $\widehat{x}$, projected gradient descent converges linearly to $\widehat{x}$ when the constraint space $\mathcal{C}$ satisfies the geometric properties described above.

- Finally, we apply these ideas to a range of specific examples: low-rank matrix estimation (where optimization is carried out under a rank constraint), sparse estimation (with nonconvex regularizers such as SCAD offering a lower-shrinkage alternative to the $\ell_1$ norm), and several other nonconvex constraints. We discuss some interesting differences between constraining versus penalizing a nonconvex regularization function, in the context of sparse estimation. For the low-rank setting, we propose an approximate projection step that provides a computationally efficient alternative for low-rank estimation problems, which we then explore empirically with simulations.

## 2   Concavity coefficients for a nonconvex constraint space

We begin by studying several properties which describe the extent to which the constraint set $\mathcal{C} \subset \mathbb{R}^d$ deviates from convexity. To quantify the concavity of $\mathcal{C}$, we will define the (global) concavity coefficient of $\mathcal{C}$, denoted $\gamma = \gamma(\mathcal{C})$, which we will later expand to local measures of concavity, $\gamma_x(\mathcal{C})$, indexed over points $x \in \mathcal{C}$. We examine several definitions of this concavity coefficient: essentially, we consider four properties that would hold if $\mathcal{C}$ were convex, and then use $\gamma$ to characterize the extent to which these properties are violated. Our definitions are closely connected to the notion of *prox-regular sets* in the analysis literature, and we will discuss this connection in detail in Section 2.3 below.

Since we are interested in developing flexible tools for high-dimensional optimization problems, several different norms will appear in the definitions of the concavity coefficients:

- The Euclidean $\ell_2$ norm, $\|\cdot\|_2$. Projections to $\mathcal{C}$ will always be taken with respect to the $\ell_2$ norm, and our later convergence guarantees will also be given with respect to this norm. If our variable is a matrix $X \in \mathbb{R}^{n \times m}$, the Euclidean $\ell_2$ norm is known as the Frobenius norm, $\|X\|_{\mathsf{F}} = \sqrt{\sum_{ij} X_{ij}^2}$.

- A "structured" norm $\|\cdot\|$, which can be chosen to be any norm on $\mathbb{R}^d$. In some cases it may be the $\ell_2$ norm, but often it will be a different norm reflecting natural structure in the problem. For instance, for a low-rank estimation problem, if $\mathcal{C}$ is a set of rank-constrained matrices then we will work with the nuclear norm, $\|\cdot\| = \|\cdot\|_{\text{nuc}}$ (defined as the sum of the singular values of the matrix). For sparse signals, we will instead use the $\ell_1$ norm, $\|\cdot\| = \|\cdot\|_1$.

- A norm $\|\cdot\|^*$, which is the dual norm to the structured norm $\|\cdot\|$. For low-rank matrix problems, if we work with the nuclear norm, $\|\cdot\| = \|\cdot\|_{\text{nuc}}$, then the dual norm is given by the spectral norm, $\|\cdot\|^* = \|\cdot\|_{\text{sp}}$ (i.e. the largest singular value of the matrix, also known as the matrix operator norm). For sparse problems, if $\|\cdot\| = \|\cdot\|_1$ then its dual is given by the $\ell_\infty$ norm, $\|\cdot\|^* = \|\cdot\|_\infty$.

When we take projections to the constraint set $\mathcal{C}$, if the minimizer $P_\mathcal{C}(z) \in \arg\min_{x \in \mathcal{C}} \|x - z\|_2$ is non-unique, then we write $P_\mathcal{C}(z)$ to denote any point chosen from this set. Throughout, any assumption or claim involving $P_\mathcal{C}(z)$ should be interpreted as holding for any choice of $P_\mathcal{C}(z)$. From this point on, we will assume without comment that $\mathcal{C}$ is closed and nonempty so that the set $\arg\min_{x \in \mathcal{C}} \|x - z\|_2$ is nonempty for any $z$.

We now present several definitions of the concavity coefficient of $\mathcal{C}$.

**Curvature** First, we define $\gamma$ as a bound on the extent to which a convex combination of two elements of $\mathcal{C}$ may lie outside of $\mathcal{C}$: for $x, y \in \mathcal{C}$,

$$\limsup_{t \searrow 0} \frac{\min_{z \in \mathcal{C}} \|z - ((1-t)x + ty)\|}{t} \leq \gamma \|x - y\|_2^2. \tag{1}$$

**Approximate contraction** Second, we define $\gamma$ via a condition requiring that the projection operator $P_\mathcal{C}$ is approximately contractive in a neighborhood of the set $\mathcal{C}$, that is, $\|P_\mathcal{C}(z) - P_\mathcal{C}(w)\|_2$ is not much larger than $\|z - w\|_2$: for $x, y \in \mathcal{C}$,

For any $z, w \in \mathbb{R}^d$ with $P_\mathcal{C}(z) = x$ and $P_\mathcal{C}(w) = y$,
$$\left(1 - \gamma\|z - x\|^* - \gamma\|w - y\|^*\right) \cdot \|x - y\|_2 \leq \|z - w\|_2. \tag{2}$$

For convenience in our theoretical analysis we will also consider a weaker "one-sided" version of this property, where one of the two points is assumed to already lie in $\mathcal{C}$: for $x, y \in \mathcal{C}$,

For any $z \in \mathbb{R}^d$ with $P_\mathcal{C}(z) = x$, $\quad \left(1 - \gamma\|z - x\|^*\right) \cdot \|x - y\|_2 \leq \|z - y\|_2. \tag{3}$

**First-order optimality** For our third characterization of the concavity coefficient, we consider the standard first-order optimality conditions for minimization over a convex set, and measure the extent to which they are violated when optimizing over $\mathcal{C}$: for $x, y \in \mathcal{C}$,[1]

For any differentiable $\mathsf{f} : \mathbb{R}^d \to \mathbb{R}$ such that $x$ is a local minimizer of $\mathsf{f}$ over $\mathcal{C}$,
$$\langle y - x, \nabla \mathsf{f}(x) \rangle \geq -\gamma \|\nabla \mathsf{f}(x)\|^* \|y - x\|_2^2. \tag{4}$$

**Inner products** Fourth, we introduce an inner product condition, requiring that projection to the constraint set $\mathcal{C}$ behaves similarly to a convex projection: for $x, y \in \mathcal{C}$,

For any $z \in \mathbb{R}^d$ with $P_\mathcal{C}(z) = x$, $\quad \langle y - x, z - x \rangle \leq \gamma \|z - x\|^* \|y - x\|_2^2. \tag{5}$

We will see later that, by choosing $\|\cdot\|$ to reflect the structure in the signal (rather than working only with the $\ell_2$ norm), we are able to obtain a more favorable scaling in our concavity coefficients, and hence to prove meaningful convergence results in high-dimensional settings. On the other hand, regardless of our choice of $\|\cdot\|$, note that the $\ell_2$ norm also appears in the definition of the concavity coefficients, as is natural when working with inner products.

Our first main result shows that the above conditions are in fact exactly equivalent:

**Theorem 1.** *The properties* (1), (2), (3), (4), *and* (5) *are equivalent; that is, for a fixed choice* $\gamma \in [0, \infty]$, *they either all hold for every* $x, y \in \mathcal{C}$, *or all fail to hold for some* $x, y \in \mathcal{C}$.

Formally, we will define $\gamma(\mathcal{C})$ to be the smallest value such that the above properties hold:

$$\gamma(\mathcal{C}) := \min\left\{\gamma \in [0, \infty] : \text{Properties (1), (2), (3), (4), (5) hold for all } x, y \in \mathcal{C}\right\}.$$

However, this global coefficient $\gamma(\mathcal{C})$ is often of limited use in practical settings, since many sets are well-behaved locally but not globally. For instance, the set $\mathcal{C} = \{X \in \mathbb{R}^{n \times m} : \text{rank}(X) \leq r\}$ has $\gamma(\mathcal{C}) = \infty$, but exhibits smooth curvature and good convergence behavior as long as we stay away from rank-degenerate matrices (that is, matrices with $\text{rank}(X) < r$). Since we may often want to ensure convergence in this type of setting where global concavity cannot be bounded, we next turn to a local version of the same concavity bounds.

---

[1] A more general form of this condition, with $\mathsf{f}$ Lipschitz but not necessarily differentiable, appears in Appendix A.2.1.

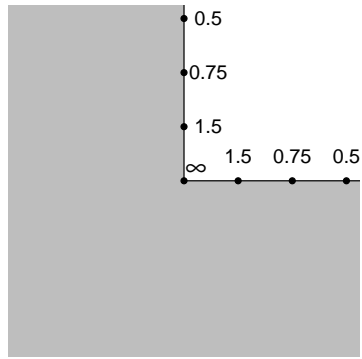

Figure 1: A simple example of the local concavity coefficients on $\mathcal{C} = \{x \in \mathbb{R}^2 : x_1 \leq 0 \text{ or } x_2 \leq 0\}$. The gray shaded area represents $\mathcal{C}$ while the numbers give the local concavity coefficients at each marked point.

## 2.1  Local concavity coefficients

We now consider the *local concavity coefficients* $\gamma_x(\mathcal{C})$, measuring the concavity in a set $\mathcal{C}$ relative to a specific point $x$ in the set. We will see examples later on where $\gamma(\mathcal{C}) = \infty$ but $\gamma_x(\mathcal{C})$ is bounded for many points $x \in \mathcal{C}$.

First we define a set of "degenerate points",
$$\mathcal{C}_{\mathsf{dgn}} = \{x \in \mathcal{C} : \ P_{\mathcal{C}} \text{ is not continuous over any neighborhood of } x\},$$
and then let
$$\gamma_x(\mathcal{C}) = \begin{cases} \infty, & x \in \mathcal{C}_{\mathsf{dgn}}, \\ \min\{\gamma \in [0,\infty] : \text{Property (*) holds for this point } x \text{ and any } y \in \mathcal{C}\}, & x \notin \mathcal{C}_{\mathsf{dgn}}, \end{cases} \quad (6)$$

where the property (*) may refer to any of the four definitions of the concavity coefficients,[2] namely (1), (3), (4), or (5). We will see shortly why it is necessary to make an exception for the degenerate points $x \in \mathcal{C}_{\mathsf{dgn}}$ in the definition of these coefficients.

Our next main result shows that the equivalence between the four properties (1), (3), (4), and (5) in terms of the global concavity coefficient $\gamma(\mathcal{C})$, holds also for the local coefficients:

**Theorem 2.** *For all $x \in \mathcal{C}$, the definition (6) of $\gamma_x(\mathcal{C})$ is equivalent for all four choices of the property (*), namely the conditions (1), (3), (4), or (5).*

To develop an intuition for the global and local concavity coefficients, we give a simple example in $\mathbb{R}^2$ (relative to the $\ell_2$ norm, i.e. $\|\cdot\| = \|\cdot\|^* = \|\cdot\|_2$), displayed in Figure 1. Define $\mathcal{C} = \{x \in \mathbb{R}^2 : x_1 \leq 0 \text{ or } x_2 \leq 0\}$. Due to the degenerate point $x = (0,0)$, we can see that $\gamma(\mathcal{C}) = \infty$ in this case. The local concavity coefficients are given by
$$\begin{cases} \gamma_x(\mathcal{C}) = \infty, & \text{if } x = (0,0), \\ \gamma_x(\mathcal{C}) = \frac{1}{2t}, & \text{if } x = (t,0) \text{ or } (0,t) \text{ for } t > 0, \\ \gamma_x(\mathcal{C}) = 0, & \text{if } x_1 < 0 \text{ or } x_2 < 0. \end{cases}$$

Note that at the degenerate point $x = (0,0)$, $\mathcal{C}$ actually contains all convex combinations of this point $x$ with any $y \in \mathcal{C}$, and so the curvature condition (1) is satisfied with $\gamma = 0$. However, $x \in \mathcal{C}_{\mathsf{dgn}}$, so we nonetheless set $\gamma_x(\mathcal{C}) = \infty$.

---

[2]In this definition, We only consider the "one-sided" formulation (3) of the contraction property, since the two-sided formulation (2) would involve the local concavity coefficient at both $x$ and $y$ due to symmetry—we will see in Lemma 4 below that a version of the two-sided contraction property still holds using local coefficients.

Practical high-dimensional examples, such as a rank constraint, will be discussed in depth in Section 5. For example we will see that, for the rank-constrained set $\mathcal{C} = \{X \in \mathbb{R}^{n \times m} : \mathrm{rank}(X) \leq r\}$, the local concavity coefficients satisfy $\gamma_X(\mathcal{C}) = \frac{1}{2\sigma_r(X)}$ relative to the nuclear norm.

In general, a rough intuition for the local coefficients is that:

- If $x$ lies in the interior of $\mathcal{C}$, or if $\mathcal{C}$ is convex, then $\gamma_x(\mathcal{C}) = 0$;

- If $x$ lies on the boundary of $\mathcal{C}$, which is a nonconvex set with a smooth boundary, then we will typically see a finite but nonzero $\gamma_x(\mathcal{C})$;

- $\gamma_x(\mathcal{C}) = \infty$ can indicate a nonconvex cusp or other degeneracy at the point $x$.

## 2.2 Properties

We next prove some properties of the local coefficients $\gamma_x(\mathcal{C})$ that will be useful for our convergence analysis, as well as for gaining intuition for these coefficients.

First, the global and local coefficients are related in the natural way:

**Lemma 1.** *For any $\mathcal{C}$, $\gamma(\mathcal{C}) = \sup_{x \in \mathcal{C}} \gamma_x(\mathcal{C})$.*

Next, observe that $x \mapsto \gamma_x(\mathcal{C})$ is not continuous in general (in particular, since $\gamma_x(\mathcal{C}) = 0$ in the interior of $\mathcal{C}$ but is often positive on the boundary). However, this map does satisfy upper semi-continuity:

**Lemma 2.** *The function $x \mapsto \gamma_x(\mathcal{C})$ is upper semi-continuous over $x \in \mathcal{C}$.*

Furthermore, setting $\gamma_x(\mathcal{C}) = \infty$ at the degenerate points $x \in \mathcal{C}_{\mathsf{dgn}}$ is natural in the following sense: the resulting map $x \mapsto \gamma_x(\mathcal{C})$ is the minimal upper semi-continuous map such that the relevant local concavity properties are satisfied. We formalize this with the following lemma:

**Lemma 3.** *For any $u \in \mathcal{C}_{\mathsf{dgn}}$, for any of the four conditions, (1), (3), (4), or (5), this property does not hold in any neighborhood of $u$ for any finite $\gamma$. That is, for any $r > 0$,*

$$\min\left\{\gamma \geq 0 : \text{Property (*) holds for all } x \in \mathcal{C} \cap \mathbb{B}_2(u, r) \text{ and for all } y \in \mathcal{C}\right\} = \infty,$$

where (*) may refer to any of the four equivalent properties, i.e. (1), (3), (4), and (5). (Here $\mathbb{B}_2(u, r)$ is the ball of radius $r$ around the point $u$, with respect to the $\ell_2$ norm.)

Finally, the next result shows that two-sided contraction property (2) holds using local coefficients, meaning that all five definitions of concavity coefficients are equivalent:

**Lemma 4.** *For any $z, w \in \mathbb{R}^d$,*

$$\left(1 - \gamma_{P_\mathcal{C}(z)}(\mathcal{C}) \|z - P_\mathcal{C}(z)\|^* - \gamma_{P_\mathcal{C}(w)}(\mathcal{C}) \|w - P_\mathcal{C}(w)\|^*\right) \cdot \|P_\mathcal{C}(z) - P_\mathcal{C}(w)\|_2 \leq \|z - w\|_2$$

In particular, for any fixed $c \in (0, 1)$, Lemma 4 proves that

$$P_\mathcal{C} \text{ is } c\text{-Lipschitz over the set } \left\{z \in \mathbb{R}^d : 2\gamma_{P_\mathcal{C}(z)}(\mathcal{C}) \|z - P_\mathcal{C}(z)\|^* \leq 1 - c\right\}, \tag{7}$$

where the Lipschitz constant is defined with respect to the $\ell_2$ norm. This provides a sort of converse to our definition of the degenerate points, where we set $\gamma_x(\mathcal{C}) = \infty$ for all $x \in \mathcal{C}_{\mathsf{dgn}}$, i.e. all points $x$ where $P_\mathcal{C}$ is *not* continuous in any neighborhood of $x$.

## 2.3 Connection to prox-regular sets

The notion of prox-regular sets and sets of positive reach arises in the literature on nonsmooth analysis in Hilbert spaces, for instance see Colombo and Thibault [10] for a comprehensive overview of the key results in this area. The work on prox-regular sets generalizes also to the notion of prox-regular functions (see e.g. Rockafellar and Wets [26, Chapter 13.F]).

A prox-regular set is a set $\mathcal{C} \subset \mathbb{R}^d$ that satisfies[3]

$$\langle y - x, z - x \rangle \leq \frac{1}{2\rho} \|z - x\|_2 \|y - x\|_2^2, \tag{8}$$

for all $x, y \in \mathcal{C}$ and all $z \in \mathbb{R}^d$ with $P_\mathcal{C}(z) = x$, for some constant $\rho > 0$. To capture the local variations in concavity over the set $\mathcal{C}$, $\mathcal{C}$ is prox-regular with respect to a continuous function $\rho : \mathcal{C} \to (0, \infty]$ if

$$\langle y - x, z - x \rangle \leq \frac{1}{2\rho(x)} \|z - x\|_2 \|y - x\|_2^2 \tag{9}$$

for all $x, y \in \mathcal{C}$ and all $z \in \mathbb{R}^d$ with $P_\mathcal{C}(z) = x$ (see e.g. Colombo and Thibault [10, Theorem 3b]).[4] Historically, prox-regularity was first formulated via the notion of "positive reach" [13]: the parameter $\rho$ appearing in (8) is the largest radius such that the projection operator $P_\mathcal{C}$ is unique for all points $z$ within distance $\rho$ of the set $\mathcal{C}$; in the local version (9), the radius is allowed to vary locally as a function of $x \in \mathcal{C}$.

These definitions (8) and (9) exactly coincide with our inner product condition (5), in the special case that $\|\cdot\|$ is the $\ell_2$ norm, by taking $\gamma = \frac{1}{2\rho}$ or, for the local coefficients, $\gamma = \frac{1}{2\rho(x)}$. In the $\ell_2$ setting, there is substantial literature exploring the equivalence between many different characterizations of prox-regularity, including properties that are equivalent to each of our characterizations of the local concavity coefficients. Here we note a few places in the literature where these conditions appear, and refer the reader to Colombo and Thibault [10] for historical background on these ideas. The curvature condition (1) is proved in Colombo and Thibault [10, Proposition 9, Theorem 14(q)]. The one- and two-sided contraction conditions (3) and (2) appear in Federer [13, Section 4.8] and Colombo and Thibault [10, Theorem 14(g))]; the inner product condition (5) can be found in Federer [13, Section 4.8], Colombo and Thibault [10, Theorem 3(b)], Canino [6, Definition 1.5], and Colombo and Marques [9, Definition 2.1]. The first-order optimality condition (4) is closely related to the inner product condition, when formulated using the ideas of normal cones and proximal normal cones (for instance, Rockafellar and Wets [26, Theorem 6.12] relates gradients of f to normal cones at $x$).

The distinctions between our definitions and results on local concavity coefficients, and the literature on prox-regularity, center on two key differences: the role of continuity, and the flexibility of the structured norm $\|\cdot\|$ (rather than the $\ell_2$ norm). We discuss these two separately.

**Continuity**  In the literature on prox-regular sets, the "reach" function $x \mapsto \rho(x) \in (0, \infty]$ is assumed to be continuous [10, Definition 1]. Equivalently, we could take a continuous function $x \mapsto \gamma_x = \frac{1}{2\rho(x)} \in [0, \infty)$ to agree with the notation of our local concavity coefficients. However, this is not the same as finding the *smallest* value $\gamma_x$ such that the concavity coefficient conditions are satisfied (locally at the point $x$). For our definitions, we do not enforce continuity of the map $x \mapsto \gamma_x$, and instead define $\gamma_x(\mathcal{C})$ as the smallest value such that the conditions are satisfied. This leads to substantial challenges in proving the equivalence of the various conditions; in Lemma 2 we prove that the map is naturally upper semi-continuous, which allows us to show the desired equivalences.

In terms of practical implications, in order to use the local concavity coefficients to describe the convergence behavior of optimization problems, it is critical that we allow for non-continuity. For instance, suppose that $\mathcal{C}$ is nonconvex, and its interior $\mathsf{Int}(\mathcal{C})$ is nonempty. For any $x \in \mathsf{Int}(\mathcal{C})$, the concavity coefficient conditions are satisfied with $\gamma_x = 0$. In particular, consider the first-order optimality condition (4): if $x \in \mathsf{Int}(\mathcal{C})$ is a local minimizer of some function $\mathsf{f}$, then $x$ is in fact the global minimizer of $\mathsf{f}(x)$ and we must have $\nabla \mathsf{f}(x) = 0$. On the other hand, since $\mathcal{C}$ is nonconvex, we must have $\gamma_x > 0$ for at least some of the points $x$ on the boundary of $\mathcal{C}$. If we do require a continuity assumption on the function $x \mapsto \gamma_x$, then we would be forced to have $\gamma_x > 0$ for some points $x \in \mathsf{Int}(\mathcal{C})$ as well. This means that $\gamma_x$ would not give a precise description of the behavior of first-order methods when constraining to $\mathcal{C}$—it would not reveal that non-global minima are impossible in the interior of the set. More generally, we will show in Lemma 5 that the local concavity coefficients (defined as the lowest possible constants, as in (6)) provide a tight characterization of the convergence behavior of projected gradient descent over the constraint set $\mathcal{C}$; if we enforce continuity, we would be forced to choose larger values for $\gamma_x(\mathcal{C})$ at some points $x \in \mathcal{C}$, and the concavity coefficients would no longer be both necessary and sufficient for convergence.

---

[3]The notion of prox-regularity is typically defined over any Hilbert space with its norm $|\cdot| = \sqrt{\langle \cdot, \cdot \rangle}$ in place of the $\ell_2$ norm, however, we restrict to the case of $\mathbb{R}^d$ for ease of comparison.

[4]There also exists in the literature an alternative notion of local prox-regularity, studied by Shapiro [27], Poliquin et al. [25], Mazade and Thibault [21], and others, where $\mathcal{C}$ is said to be prox-regular at a point $u \in \mathcal{C}$ if (8) holds over $x, y \in \mathcal{C}$ that are in some arbitrarily small neighborhood of $u$. Importantly, this definition of local prox-regularity differs from our notion of local concavity coefficients, since we allow $y$ to range over the entire set; this distinction is critical for studying convergence to a *global* rather than local minimizer.

One related point is that, by allowing for $\gamma_x(\mathcal{C})$ to be infinite if needed (which would be equivalent to allowing the "reach" $\rho(x)$ to be zero for some $x$), we can accommodate constraint sets such as the low-rank matrix constraint, $\mathcal{C} = \{X \in \mathbb{R}^{n \times m} : \mathrm{rank}(X) \leq r\}$. Recalling that $\gamma_X(\mathcal{C}) = \frac{1}{2\sigma_r(X)}$ as mentioned earlier, we see that a rank-deficient matrix $X$ (i.e. $\mathrm{rank}(X) < r$) will have $\gamma_X(\mathcal{C}) = \infty$. By not requiring that the concavity coefficient is finite (equivalently, that the reach is positive), we avoid the need for any inelegant modifications (e.g. working with a truncated set such as $\mathcal{C} = \{X : \mathrm{rank}(X) \leq r, \sigma_r(X) \geq \epsilon\}$).

**Structured norms**  Prox-regularity (or equivalently the notion of positive reach) is studied in the literature in a Hilbert space, with respect to its norm, which in $\mathbb{R}^d$ means the $\ell_2$ norm (or a weighted $\ell_2$ norm).[5] In contrast, our work defines local concavity coefficients with respect to a general structured norm $\|\cdot\|$, such as the $\ell_1$ norm in a sparse signal estimation setting. To see the distinction, compare our inner product condition (5) with the definition of prox-regularity (8).

Of course, the equivalence of all norms on $\mathbb{R}^d$ means that if $\gamma(\mathcal{C})$ is finite when defined with respect to the $\ell_2$ norm (i.e. $\mathcal{C}$ is prox-regular), then it is finite with respect to any other norm—so the importance of the distinction may not be immediately clear. As an example, let $\gamma^{\ell_1}(\mathcal{C})$ and $\gamma^{\ell_2}(\mathcal{C})$ denote the concavity coefficients with respect to the $\ell_1$ and $\ell_2$ norms. Since $\|\cdot\|_2 \leq \|\cdot\|_1 \leq \sqrt{d}\|\cdot\|_2$, we could trivially show that

$$\gamma^{\ell_2}(\mathcal{C}) \leq \gamma^{\ell_1}(\mathcal{C}) \leq \sqrt{d} \cdot \gamma^{\ell_2}(\mathcal{C}),$$

but the factor $\sqrt{d}$ is unfavorable so in many settings this is a very poor bound on $\gamma^{\ell_1}(\mathcal{C})$.

We may then ask, why can we not simply define the coefficients in terms of the $\ell_2$ norm? The reason is that in optimization problems arising in high-dimensional settings (for instance, high-dimensional regression in statistics), structured norms such as the $\ell_1$ norm (for problems involving sparse signals) or the nuclear norm (for low-rank signals) allow for statistical and computational analyses that would not be possible with the $\ell_2$ norm. In particular, we will see later on that convergence for the minimization problem $\min_{x \in \mathcal{C}} \mathsf{g}(x)$ will depend on bounding $\|\nabla \mathsf{g}(x)\|^*$. If $\|\cdot\|$ is the $\ell_1$ norm, for instance, then $\|\nabla \mathsf{g}(x)\|^* = \|\nabla \mathsf{g}(x)\|_\infty$ will in general be much smaller than $\|\nabla \mathsf{g}(x)\|_2$. For instance, in a statistical problem, if $\nabla \mathsf{g}(x)$ consists of Gaussian or sub-gaussian noise at the true parameter vector $x$, then $\|\nabla \mathsf{g}(x)\|_\infty \sim \sqrt{\log(d)}$ while $\|\nabla \mathsf{g}(x)\|_2 \sim \sqrt{d}$. Therefore, being able to bound the concavity of $\mathcal{C}$ with respect to the $\ell_1$ norm rather than the $\ell_2$ norm is crucial for analyzing convergence in a high-dimensional setting.

In the next section, we will study how the choice of the norm $\|\cdot\|$ and its dual $\|\cdot\|^*$ relates to the convergence properties of projected gradient descent.

---

[5] The notion of prox-regular sets has been extended to Banach spaces, e.g. by Bernard et al. [1], but this does not relate directly to our results.

# 3   Fast convergence of projected gradient descent

Consider an optimization problem constrained to a nonconvex set, $\min\{\mathsf{g}(x) : x \in \mathcal{C}\}$, where $\mathsf{g} : \mathbb{R}^d \to \mathbb{R}$ is a differentiable function. We will work with projected gradient descent algorithms in the setting where $\mathsf{g}$ is convex or approximately convex, while $\mathcal{C}$ is nonconvex with local concavity coefficients $\gamma_x(\mathcal{C})$. After choosing some initial point $x_0 \in \mathcal{C}$, for each $t \geq 0$ we define

$$\begin{cases} x'_{t+1} = x_t - \eta \nabla \mathsf{g}(x_t), \\ x_{t+1} = P_\mathcal{C}(x'_{t+1}), \end{cases} \tag{10}$$

where if $P_\mathcal{C}(x'_{t+1})$ is not unique then any closest point may be chosen.

## 3.1   Assumptions

**Assumptions on $\mathsf{g}$**    We first consider the objective function $\mathsf{g}$. Let $\widehat{x}$ be the target of our optimization procedure, $\widehat{x} \in \arg\min_{x \in \mathcal{C}} \mathsf{g}(x)$.

We assume that $\mathsf{g}$ satisfies restricted strong convexity (RSC) and restricted smoothness (RSM) conditions over $x, y \in \mathcal{C}$,

$$\mathsf{g}(y) \geq \mathsf{g}(x) + \langle y - x, \nabla \mathsf{g}(x) \rangle + \frac{\alpha}{2} \|x - y\|_2^2 - \frac{\alpha}{2} \varepsilon_\mathsf{g}^2, \tag{11}$$

and

$$\mathsf{g}(y) \leq \mathsf{g}(x) + \langle y - x, \nabla \mathsf{g}(x) \rangle + \frac{\beta}{2} \|x - y\|_2^2 + \frac{\alpha}{2} \varepsilon_\mathsf{g}^2. \tag{12}$$

Without loss of generality we can take $\alpha \leq \beta$. As is common in the low-rank factorized optimization literature, we will work in a local neighborhood of the target $\widehat{x}$ by assuming that our initialization point lies within radius $\rho$ of $\widehat{x}$, which will allow us to require these conditions on $\mathsf{g}$ to hold only locally.

The term $\varepsilon_\mathsf{g}$ gives some "slack" in our assumption on $\mathsf{g}$, and is intended to capture some vanishingly small error level. This term is often referred to as the "statistical error" in the high-dimensional statistics literature, which represents the best-case scaling of the accuracy of our recovered solution. Often $\widehat{x}$ may represent a global minimizer which is within radius $\varepsilon_\mathsf{g}$ of some "true" parameter in a statistical setting; therefore, converging to $\widehat{x}$ up to an error of magnitude $\varepsilon_\mathsf{g}$ means that the recovered solution is as accurate as $\widehat{x}$ at recovering the true parameter. For instance, often we will have $\varepsilon_\mathsf{g} \sim \sqrt{\frac{\log(d)}{n}}$ in a statistical setting where we are solving a sparse estimation problem of dimension $d$ with sample size $n$.

**Assumptions on $\mathcal{C}$**    Next, turning to the nonconvexity of $\mathcal{C}$, we will assume local concavity coefficients $\gamma_x(\mathcal{C})$ that are not too large in a neighborhood of $\widehat{x}$, with details given below. We furthermore assume a *norm compatibility condition*,

$$\|z - P_\mathcal{C}(z)\|^* \leq \phi \min_{x \in \mathcal{C}} \|z - x\|^* \text{ for all } z \in \mathbb{R}^d, \tag{13}$$

for some constant $\phi \geq 1$. The norm compatibility condition is trivially true with $\phi = 1$ if $\|\cdot\|$ is the $\ell_2$ norm, since $P_\mathcal{C}$ is a projection with respect to the $\ell_2$ norm. We will see that in many natural settings it holds even for other norms, often with $\phi = 1$.

**Assumptions on gradient and initialization**    Finally, we assume a gradient condition that reveals the connection between the curvature of the nonconvex set $\mathcal{C}$ and the target function $\mathsf{g}$: we require that

$$2\phi \cdot \max_{x, x' \in \mathcal{C} \cap \mathbb{B}_2(\widehat{x}, \rho)} \gamma_x(\mathcal{C}) \|\nabla \mathsf{g}(x')\|^* \leq (1 - c_0) \cdot \alpha. \tag{14}$$

(Since $x \mapsto \gamma_x(\mathcal{C})$ is upper semi-continuous, if $\mathsf{g}$ is continuously differentiable, then we can find some radius $\rho > 0$ and some constant $c_0 > 0$ satisfying this condition, as long as $2\phi \gamma_{\widehat{x}}(\mathcal{C}) \|\nabla \mathsf{g}(\widehat{x})\|^* < \alpha$.) Our projected gradient descent algorithm will then succeed if initialized within this radius $\rho$ from the target point $\widehat{x}$, with an appropriate step size. We will discuss the necessity of this type of initialization condition below in Section 3.4.

In practice, relaxing the constraint $x \in \mathcal{C}$ to a convex constraint (or convex penalty) is often sufficient for providing a good initialization point. For example, in low-rank matrix setting, if we would like to solve $\arg\min\{\mathsf{g}(X) : \mathrm{rank}(X) \leq r\}$, we may first solve $\arg\min_X\{\mathsf{g}(X) + \lambda \|X\|_{\mathrm{nuc}}\}$, where $\|X\|_{\mathrm{nuc}}$ is the nuclear norm and $\lambda \geq 0$ is a penalty parameter (which we would tune to obtain the desired rank for $X$). Alternately, in some settings, it may be sufficient to solve an unconstrained problem $\arg\min_X \mathsf{g}(X)$ and then project to the constraint set, $P_\mathcal{C}(X)$. For some detailed examples of suitable initialization procedures for various low-rank matrix estimation problems, see e.g. Chen and Wainwright [8], Tu et al. [30].

## 3.2 Convergence guarantee

We now state our main result, which proves that under these conditions, initializing at some $x_0 \in \mathcal{C}$ sufficiently close to $\widehat{x}$ will guarantee fast convergence to $\widehat{x}$.

**Theorem 3.** *Let $\mathcal{C} \subset \mathbb{R}^d$ be a constraint set and let $\mathsf{g}$ be a differentiable function, with minimizer $\widehat{x} \in \arg\min_{x \in \mathcal{C}} \mathsf{g}(x)$. Suppose $\mathcal{C}$ satisfies the norm compatibility condition* (13) *with parameter $\phi$, and $\mathsf{g}$ satisfies restricted strong convexity* (11) *and restricted smoothness* (12) *with parameters $\alpha, \beta, \varepsilon_\mathsf{g}$ for all $x, y \in \mathbb{B}_2(\widehat{x}, \rho)$, and the initialization condition* (14) *for some $c_0 > 0$. If the initial point $x_0 \in \mathcal{C}$ and the error level $\varepsilon_\mathsf{g}$ satisfy $\|x_0 - \widehat{x}\|_2^2 < \rho^2$ and $\varepsilon_\mathsf{g}^2 < \frac{c_0 \rho^2}{1.5}$, then for each step $t \geq 0$ of the projected gradient descent algorithm* (10) *with step size $\eta = 1/\beta$,*

$$\|x_t - \widehat{x}\|_2^2 \leq \left(1 - c_0 \cdot \frac{2\alpha}{\alpha + \beta}\right)^t \|x_0 - \widehat{x}\|_2^2 + \frac{1.5\varepsilon_\mathsf{g}^2}{c_0}.$$

In other words, the iterates $x_t$ converge linearly to the minimizer $\widehat{x}$, up to precision level $\varepsilon_\mathsf{g}$.

## 3.3 Comparison to related work

We now compare to several related results for convex and nonconvex projected gradient descent. (For methods that are specific to the problem of optimization over low-rank matrices, we will discuss this comparison and perform simulations later on.)

**Comparison to convex optimization**    To compare this result to the convex setting, if $\mathcal{C}$ is a convex set and $\mathsf{g}$ is $\alpha$-strongly convex and $\beta$-smooth, then we can set $c_0 = 1$ and $\varepsilon_\mathsf{g} = 0$. Our result then yields

$$\|x_t - \widehat{x}\|_2^2 \leq \left(1 - \frac{2\alpha}{\alpha + \beta}\right)^t \|x_0 - \widehat{x}\|_2^2 = \left(\frac{\beta - \alpha}{\beta + \alpha}\right)^t \|x_0 - \widehat{x}\|_2^2,$$

matching known rates for the convex setting (see e.g. Bubeck [2, Theorem 3.10]).

**Comparison to known results using descent cones**    Oymak et al. [23] study projected gradient descent for a linear regression setting, $\mathsf{g}(x) = \frac{1}{2}\|b - Ax\|_2^2$, while constraining some potentially nonconvex regularizer, $\mathcal{C} = \{x : \mathrm{Pen}(x) \leq c\}$. Given a true solution $x^\star \in \mathcal{C}$ (for instance, in a statistical setting, we may have $b = Ax^\star + \text{(noise)}$), their work focuses on the *descent cone* of $\mathcal{C}$ at $x^\star$, given by

$$\mathrm{DC}_{x^\star} = \text{Smallest closed cone containing } \{u : \mathrm{Pen}(x^\star + u) \leq c\}.$$

(Trivially we will have $x_t - x^\star \in \mathrm{DC}_{x^\star}$ since $x_t \in \mathcal{C}$.) Their results characterize the convergence of projected gradient descent in terms of the eigenvalues of $A^\top A$ restricted to this cone. For simplicity, we show their result specialized to the noiseless setting, i.e. $b = Ax^\star$, given in [23, Theorem 1.2]:

$$\|x_t - x^\star\|_2 \leq \left(2 \cdot \max_{u,v \in \mathrm{DC}_{x^\star} \cap \mathbb{S}^{d-1}} u^\top \left(\mathbf{I}_d - \eta A^\top A\right) v\right)^t \|x^\star\|_2. \tag{15}$$

For this result to be meaningful we of course need the radius of convergence to be $< 1$. For a convex constraint set $\mathcal{C}$ (i.e. if $\text{Pen}(x)$ is convex), the factor of 2 can be removed. In the nonconvex setting, however, the factor of 2 means that the maximum in (15) must be $< \frac{1}{2}$ for the bound to ensure convergence. Noting that $\nabla \mathsf{g}(x) = A^\top A(x - x^\star)$ in this problem, by setting $u = v \propto x - x^\star$ we see that (15) effectively requires that $\eta > \frac{1}{2\alpha}$, where $\alpha$ is the restricted strong convexity parameter (11). However, we also know that $\eta \leq \frac{1}{\beta}$ is generally a necessary condition to ensure stability of projected gradient descent; if $\eta > \frac{1}{\beta}$ then we may see values of $\mathsf{g}$ increase over iterations, i.e. $\mathsf{g}(x_1) > \mathsf{g}(x_0)$. Therefore, the condition (15) effectively requires that $\mathsf{g}$ is well-conditioned with $\beta \lesssim 2\alpha$, and furthermore that $x^\star$ is not in the interior of $\mathcal{C}$ (since, if this were the case, then $\text{DC}_{x^\star} = \mathbb{R}^d$). On the other hand, if the radius in (15) is indeed $< 1$, then their work does not assume any type of initialization condition for convergence to be successful, in contrast to our initialization assumption (14).

**Comparison to known results for iterative hard thresholding**  We now compare our results to those of Jain et al. [17], which specifically treat the iterative hard thresholding algorithm for a sparsity constraint or a rank constraint,
$$\mathcal{C} = \{x \in \mathbb{R}^d : |\text{support}(x)| \leq k\} \text{ or } \mathcal{C} = \{X \in \mathbb{R}^{n \times m} : \text{rank}(X) \leq r\}.$$
In their work, they take a substantially different approach: instead of bounding the distance between $x_t$ and the minimizer $\widehat{x} \in \arg\min_{x \in \mathcal{C}} \mathsf{g}(x)$, they instead take $\widehat{x}$ to be a minimizer over a stronger constraint,
$$\widehat{x} = \arg\min_{|\text{support}(x)| \leq k^\star} \mathsf{g}(x) \quad \text{or} \quad \widehat{X} = \arg\min_{\text{rank}(X) \leq r^\star} \mathsf{g}(X),$$
taking $k^\star \ll k$ or $r^\star \ll r$ to enforce that the sparsity of $\widehat{x}$ or rank of $\widehat{X}$ is much lower than the optimization constraint set $\mathcal{C}$. With this definition, then bound the gap in objective function values, $\mathsf{g}(x_t) - \mathsf{g}(\widehat{x})$. In other words, the objective function value $\mathsf{g}(x_t)$ is, up to a small error, no larger than the best value obtained over the substantially more restricted set of $k^\star$-sparse vectors or of rank-$r^\star$ matrices. By careful use of this gap $k^\star \ll k$ or $r^\star \ll r$, their analysis allows for convergence results from *any* initialization point $x_0 \in \mathcal{C}$. In contrast, our work allows $\widehat{x}$ to lie anywhere in $\mathcal{C}$, but this comes at the cost of assuming a *local* initialization point $x_0 \in \mathcal{C} \cap \mathbb{B}_2(\widehat{x}, \rho)$.

This result suggests a possible two-phase approach: first, we might optimize over a larger rank constraint $\mathcal{C} = \{X : \text{rank}(X) \leq k\}$ where $k \gg k^\star$ to obtain the convergence guarantees of Jain et al. [17] (which do not assume a good initialization point, but obtain weaker guarantees); then, given the solution over rank $k$ as a good initialization point, we would then optimize over the tighter constraint $\mathcal{C} = \{X : \text{rank}(X) \leq k^\star\}$ to obtain our stronger guarantees.

**Comparison to results on prox-regular functions**  Pennanen [24] studies conditions for linear convergence of the *proximal point method* for minimizing a function $\mathsf{f}(x)$, and shows that prox-regularity of $\mathsf{f}(x)$ is sufficient; Lewis and Wright [19] also study this problem in a more general setting. For our optimization problem, this translates to setting $\mathsf{f}(x) = \mathsf{g}(x) + \delta_{\mathcal{C}}(x)$, where
$$\delta_{\mathcal{C}}(x) = \begin{cases} 0, & x \in \mathcal{C}, \\ \infty, & x \notin \mathcal{C}. \end{cases}$$

(This is usually called the "indicator function" for the set $\mathcal{C}$.) If $\mathsf{g}(x) + \frac{\mu}{2}\|x\|_2^2$ is convex (i.e. the concavity of $\mathsf{g}$ is bounded) and $\mathcal{C}$ is a prox-regular set (i.e. $\gamma(\mathcal{C}) < \infty$, see Section 2.3), then $\mathsf{f}(x)$ is a prox-regular function. This work was extended by Iusem et al. [15] and others to an inexact proximal point method, allowing for error in each interation, which can be formulated to encompass the projected gradient descent algorithm studied here. Our first convergence result Theorem 3 extends these results into a high-dimensional setting by using the structured norm $\|\cdot\|$ and its dual $\|\cdot\|^*$ (e.g. the $\ell_1$ norm and its dual the $\ell_\infty$ norm), and requiring only restricted strong convexity and restricted smoothness on $\mathsf{g}$, without which we would not be able to obtain convergence guarantees in settings such as high-dimensional sparse regression or low-rank matrix estimation.

### 3.4   Initialization point and the gradient assumption

In this result, we assume that the initialization point $x_0$ is within some radius $\rho$ of the target $\widehat{x}$, ensuring that $2\phi\gamma_x(\mathcal{C})\|\nabla\mathsf{g}(x)\|^* < \alpha$ for all $x$ in the initialization neighborhood, where $\alpha$ is the restricted strong convexity (11) parameter. This type of assumption arises in much of the related literature; for example in the setting of optimization over low-rank matrices, as we will see in Section 5.1, we will require that $\|X_0 - \widehat{X}\|_\mathsf{F} \lesssim \sigma_r(\widehat{X})$, which is the same condition found in existing work such as Chen and Wainwright [8].

In fact, the following result demonstrates that the bound (14) is in a sense necessary:

**Lemma 5.** *For any constraint set $\mathcal{C}$ and any point $x \in \mathcal{C}\backslash\mathcal{C}_{\mathsf{dgn}}$ with $\gamma_x(\mathcal{C}) > 0$, for any $\alpha, \epsilon > 0$ there exists an $\alpha$-strongly convex $\mathsf{g}$ such that*

- *The gradient condition (14) is nearly satisfied at $x$, with $2\gamma_x(\mathcal{C})\|\mathsf{g}(x)\|^* \leq \alpha(1+\epsilon)$,*
- *And, $x$ is a stationary point of the projected gradient descent algorithm (10) for all sufficiently small step sizes $\eta > 0$,*
- *But $x$ does not minimize $\mathsf{g}$ over $\mathcal{C}$.*

That is, if projected gradient descent is initialized at the point $x$, then the algorithm will never leave this point, even though it is not optimal (i.e. $x$ is not the global minimizer).

We can see with a concrete example that the condition (14) may be even more critical than this lemma suggests: without this bound, we may find that projected gradient descent becomes trapped at a stationary point $x$ which is not even a *local* minimum, as in the following example.

**Example 1.** *Let $\mathcal{C} = \{X \in \mathbb{R}^{2\times 2} : \mathrm{rank}(X) \leq 1\}$, let $\mathsf{g}(X) = \frac{1}{2}\left\|X - \begin{pmatrix} 1 & 0 \\ 0 & 1+\epsilon \end{pmatrix}\right\|_\mathsf{F}^2$, and let $X_0 = \begin{pmatrix} 1 & 0 \\ 0 & 0 \end{pmatrix}$. Then trivially, we can see that $\mathsf{g}$ is $\alpha$-strongly convex for $\alpha = 1$, and that $X_0$ is a stationary point of the projected gradient descent algorithm (10) for any step size $\eta < \frac{1}{1+\epsilon}$. However, for any $0 < t < \sqrt{2\epsilon}$, setting $X = \begin{pmatrix} 1 & t \\ t & t^2 \end{pmatrix} \in \mathcal{C}$, we can see that $\mathsf{g}(X) < \mathsf{g}(X_0)$—that is, $X_0$ is stationary point, but is not a local minimum.*

*We will later calculate that $\gamma_{X_0}(\mathcal{C}) = \frac{1}{2\sigma_1(X_0)} = \frac{1}{2}$ relative to the nuclear norm $\|\cdot\| = \|\cdot\|_\mathsf{nuc}$, with norm compatibility constant $\phi = 1$ (see Section 5.1 for this calculation). Comparing against the condition (14) on the gradient of $\mathsf{g}$, since the dual norm to $\|\cdot\|_\mathsf{nuc}$ is the matrix spectral norm $\|\cdot\|_\mathsf{sp}$, we see that*

$$2\phi\gamma_{X_0}(\mathcal{C}) \cdot \|\nabla\mathsf{g}(X_0)\|_\mathsf{sp} = 2 \cdot 1 \cdot \frac{1}{2} \cdot \left\|-\begin{pmatrix} 0 & 0 \\ 0 & 1+\epsilon \end{pmatrix}\right\|_\mathsf{sp} = 1 + \epsilon = \alpha(1+\epsilon).$$

*Therefore, when the initial gradient condition (14) is even slightly violated in this example (i.e. small $\epsilon > 0$), the projected gradient descent algorithm can become trapped at a point that is not even a local minimum.*

While we might observe that in this particular example, the "bad" stationary point $X_0$ could be avoided by increasing the step size, in other settings if $\mathsf{g}$ has strong curvature in some directions (i.e. the smoothness parameter $\beta$ is large), then we cannot afford a large step size $\eta$ as it can cause the algorithm to fail to converge.

## 4   Convergence analysis using approximate projections

In some settings, computing projections $P_\mathcal{C}(x'_{t+1})$ at each step of the projected gradient descent algorithm may be prohibitively expensive; for instance in a low-rank matrix optimization problem of dimension $d \times d$, this would generally involve taking the singular value decomposition of a dense $d \times d$ matrix at each step. In these cases we may sometimes have access to a fast but approximate computation of this projection, which may come at the cost of slower convergence.

We now generalize to the idea of a *family of approximate projections*, which allows for operators that approximate projection to $\mathcal{C}$. Specifically, the approximations are carried out locally:

$$\begin{cases} x'_{t+1} = x_t - \eta \nabla \mathsf{g}(x_t), \\ x_{t+1} = P_{x_t}(x'_{t+1}), \end{cases} \tag{16}$$

where $P_{x_t}$ comes from a family of operators $P_x : \mathbb{R}^d \to \mathcal{C}$ indexed by $x \in \mathcal{C}$. Intuitively, we think of $P_x(z)$ as providing a very accurate approximation to $P_\mathcal{C}(z)$ locally for $z$ near $x$, but it may distort the projection more as we move farther away.

To allow for our convergence analysis to carry through even with these approximate projections, we assume that the family of operators $\{P_x\}$ satisfies a relaxed inner product condition:

For any $x \in \mathcal{C}$ and $z \in \mathbb{R}^d$ with $x, P_x(z) \in \mathbb{B}_2(\widehat{x}, \rho)$,

$$\langle \widehat{x} - P_x(z), z - P_x(z) \rangle \leq \max\{\underbrace{\|z - P_x(z)\|^*}_{\text{concavity term}}, \underbrace{\|z - x\|^*}_{\text{distortion term}}\} \cdot \left( \underbrace{\gamma^{\mathrm{c}} \|\widehat{x} - P_x(z)\|_2^2}_{\text{concavity term}} + \underbrace{\gamma^{\mathrm{d}} \|\widehat{x} - x\|_2^2}_{\text{distortion term}} \right). \tag{17}$$

Here the "concavity" terms are analogous to the inner product bound in (5) for exact projection to the nonconvex set $\mathcal{C}$, except with the projection $P_\mathcal{C}$ replaced by the operator $P_x$; the "distortion" terms mean that as we move farther away from $x$ the bound becomes looser, as $P_x$ becomes a less accurate approximation to $P_\mathcal{C}$.

We now present a convergence guarantee nearly identical to the result for the exact projection case, Theorem 3. We first need to state a version of the norm compatibility condition, modified for approximate projections:

$$\|z - P_x(z)\|^* \leq \phi \|z - x\|^* \text{ for all } x \in \mathcal{C} \cap \mathbb{B}_2(\widehat{x}, \rho) \text{ and } z \in \mathbb{R}^d. \tag{18}$$

We also require a modified initialization condition,

$$2\phi(\gamma^{\mathrm{c}} + \gamma^{\mathrm{d}}) \max_{x \in \mathcal{C} \cap \mathbb{B}_2(\widehat{x}, \rho)} \|\nabla \mathsf{g}(x)\|^* \leq (1 - c_0)\alpha, \tag{19}$$

and a modified version of local uniform continuity (compare to (7) for exact projections),

For any $x \in \mathcal{C} \cap \mathbb{B}_2(\widehat{x}, \rho)$, for any $\epsilon > 0$, there exists a $\delta > 0$ such that,

for any $z, w \in \mathbb{R}^d$ such that $P_x(z) \in \mathbb{B}_2(\widehat{x}, \rho)$ and $2(\gamma^{\mathrm{c}} + \gamma^{\mathrm{d}}) \|z - P_x(z)\|^* \leq 1 - c_0$,

if $\|z - w\|_2 \leq \delta$ then $\|P_x(z) - P_x(w)\|_2 \leq \epsilon$. $\tag{20}$

Our result for this setting now follows.

**Theorem 4.** *Let $\mathcal{C} \subset \mathbb{R}^d$ be a constraint set and let $\mathsf{g}$ be a differentiable function, with minimizer $\widehat{x} \in \arg\min_{x \in \mathcal{C}} \mathsf{g}(x)$. Let $\{P_x\}$ be a family of operators satisfying the inner product condition (17), the norm compatibility condition (18), and the local continuity condition (20) with parameters $\gamma^{\mathrm{c}}, \gamma^{\mathrm{d}}, \phi$ and radius $\rho$. Assume that $\mathsf{g}$ satisfies restricted strong convexity (11) and restricted smoothness (12) with parameters $\alpha, \beta, \varepsilon_\mathsf{g}$ for all $x, y \in \mathbb{B}_2(\widehat{x}, \rho)$, and the initialization condition (19) for some $c_0 > 0$. If the initial point $x_0 \in \mathcal{C}$ and the error level $\varepsilon_\mathsf{g}$ satisfy $\|x_0 - \widehat{x}\|_2^2 < \rho^2$ and $\varepsilon_\mathsf{g}^2 < \frac{c_0 \rho^2}{1.5}$, then for each step $t \geq 0$ of the approximate projected gradient descent algorithm (16) with step size $\eta = 1/\beta$,*

$$\|x_t - \widehat{x}\|_2^2 \leq \left(1 - c_0 \cdot \frac{2\alpha}{\alpha + \beta}\right)^t \|x_0 - \widehat{x}\|_2^2 + \frac{1.5\varepsilon_\mathsf{g}^2}{c_0}.$$

This convergence rate is identical to that obtained in Theorem 3 for exact projections—the only differences lie in the assumptions.

## 4.1 Exact versus approximate projections

To compare the two settings we have considered, exact projections $P_\mathcal{C}$ versus approximate projections $P_x$, we focus on a local form of the inner product condition (17) for the family of approximate operators $\{P_x\}$, rewritten to be analogous to the inner product condition (5) for exact projections. Suppose that $\gamma_u^c(\mathcal{C})$ and $\gamma_u^d(\mathcal{C})$ satisfy the property that

For any $x, y \in \mathcal{C}$ and any $z \in \mathbb{R}^d$, writing $u = P_x(z)$,
$$\langle y - u, z - u \rangle \leq \max\{\underbrace{\|z-u\|^*}_{\text{concavity term}}, \underbrace{\|z-x\|^*}_{\text{distortion term}}\} \cdot \left(\underbrace{\gamma_u^c(\mathcal{C})\|y-u\|_2^2}_{\text{concavity term}} + \underbrace{\gamma_u^d(\mathcal{C})\|y-x\|_2^2}_{\text{distortion term}}\right), \quad (21)$$

where $u \mapsto \gamma_u^c(\mathcal{C})$ and $u \mapsto \gamma_u^d(\mathcal{C})$ are upper semi-continuous maps. We now prove that the existence of a family of operators $\{P_x\}$ satisfying this general condition (21) is in fact equivalent to bounding the local concavity coefficients of $\mathcal{C}$.

**Lemma 6.** *Consider a constraint set $\mathcal{C} \subset \mathbb{R}^d$ and a norm $\|\cdot\|$ on $\mathbb{R}^d$ with dual norm $\|\cdot\|^*$. If $\mathcal{C}$ has local concavity coefficients given by $\gamma_x(\mathcal{C})$ for all $x \in \mathcal{C}$, then by defining operators $P_x = P_\mathcal{C}$ for all $x \in \mathcal{C}$, the inner product condition (21) holds with $\gamma_x^c(\mathcal{C}) = \gamma_x(\mathcal{C})$ and $\gamma_x^d(\mathcal{C}) = 0$. Conversely, if there is some family of operators $\{P_x\}_{x \in \mathcal{C}}$ satisfying the inner product condition (21), then the local concavity coefficients of $\mathcal{C}$ satisfy $\gamma_x(\mathcal{C}) \leq \gamma_x^c(\mathcal{C}) + \gamma_x^d(\mathcal{C})$, provided that $x \mapsto \gamma_x^c(\mathcal{C}), x \mapsto \gamma_x^d(\mathcal{C})$ are upper semi-continuous, and that $P_x$ also satisfies a local continuity assumption,*

$$\text{If } \gamma_x^c(\mathcal{C}) + \gamma_x^d(\mathcal{C}) < \infty \text{ and } z_t \to x, \text{ then } P_x(z_t) \to x. \quad (22)$$

For this reason, we see that generalizing from exact projection $P_\mathcal{C}$ to a family of operators $\{P_x\}$ does not expand the class of problems whose convergence is ensured by our theory; essentially, if using the approximate projection operators $P_x$ guarantees fast convergence, then the same would also be true using exact projection $P_\mathcal{C}$. However, there may be substantial computational gain in switching from exact to approximate projection, which comes with little or no cost in terms of convergence guarantees.

## 5 Examples

In this section we consider a range of nonconvex constraints arising naturally in high-dimensional statistics, and show that these sets come equipped with well-behaved local concavity coefficients (thus allowing for fast convergence of gradient descent, for appropriate functions g).

## 5.1 Low rank optimization

Estimating a matrix with low rank structure arises in a variety of problems in high-dimensional statistics and machine learning. A partial list includes PCA (principal component analysis), factor models, matrix completion, and reduced rank regression. The past few years have seen extensive results on the specific problem of optimization over the space of low-rank matrices:

$$\min\{g(X) : X \in \mathbb{R}^{n \times m}, \text{rank}(X) \leq r\},$$

where in various settings g($X$) may represent a least-squares loss from a linear matrix sensing problem, an objective function for the matrix completion problem, or a more general function satisfying some type of restricted convexity assumption. In addition to extensive earlier work on convex relaxations of this problem via the nuclear norm and other penalties, more recently this problem has been studied using the exact rank-$r$ constraint. The recent literature has generally treated the rank-constrained problem in one of two ways. First, the *iterative hard thresholding* method (also discussed earlier in Section 3.3) proceeds by taking gradient descent in the space of $n \times m$ matrices, then at each step projecting to the nearest rank-$r$ matrix in order to enforce a rank constraint on $X$. This amounts to optimizing the function g$(X)$ over the nonconvex constraint space of rank-$r$ matrices. Convergence results for this setting have been proved by Meka et al. [22], Jain et al. [17]. However, in high dimensions, a computational drawback of this method is the need to take a singular value decomposition of a (potentially dense) $n \times m$ matrix at each step. Alternately, one can consider the *factorized* approach, which reparametrizes the problem by taking a low-rank factorization, $X = AB^\top$ where $A \in \mathbb{R}^{n \times r}$ and $B \in \mathbb{R}^{m \times r}$, and pursuing alternating minimization or alternating gradient descent on the factors $A$ and $B$. Recent results in this line of work include Gunasekar et al. [14], Jain et al. [16], Sun and Luo [28], Tu et al. [30], Zhao et al. [32], Zhu et al. [34], and many others. This reformulation of the problem now consists of a highly nonconvex objective function g$(AB^\top)$ optimized over a generally convex space of factor matrices $(A, B) \in \mathbb{R}^{n \times r} \times \mathbb{R}^{m \times r}$, via alternating gradient descent or alternating minimization over the factors $A$ and $B$. In the special case where $X$ is positive semidefinite, we can instead optimize g$(AA^\top)$ via gradient descent on $A \in \mathbb{R}^{n \times r}$, which is again a nonconvex objective function being minimized over a convex space, and has also been extensively studied, e.g. by Zheng and Lafferty [33], Candès et al. [5], Chen and Wainwright [8], among others. For both of these cases, the analysis of the optimization problem is complicated by the issue of identifiability, where the factor(s) can only be identified up to rotation.

In this section, we will study the set of rank-constrained matrices
$$\mathcal{C} = \{X \in \mathbb{R}^{n \times m} : \mathrm{rank}(X) \leq r\}$$
to determine how our general framework of local concavity applies to this specific low rank setting. To avoid triviality, we assume $r < \min\{n, m\}$. Writing $\sigma_1(X) \geq \sigma_2(X) \geq \ldots$ to denote the sorted singular values of any matrix $X$, we compute the curvature condition and norm compatibility condition of $\mathcal{C}$:

**Lemma 7.** *Let $\mathcal{C} = \{X \in \mathbb{R}^{n \times m} : \mathrm{rank}(X) \leq r\}$. Then $\mathcal{C}$ has local concavity coefficients given by $\gamma_X(\mathcal{C}) = \frac{1}{2\sigma_r(X)}$ for all $X \in \mathcal{C}$, and satisfies the norm compatibility condition* (13) *with parameter $\phi = 1$, with respect to norms $\|\cdot\| = \|\cdot\|_{\mathrm{nuc}}$ and $\|\cdot\|^* = \|\cdot\|_{\mathrm{sp}}$.*

Thus, as long as the objective function g satisfies the appropriate conditions, we can expect projected gradient descent over the space of rank-$r$ matrices to converge well when we initialize at some matrix $X_0$ that is within a distance smaller than $\sigma_r(\widehat{X})$ from the target matrix $\widehat{X}$, so that $\gamma_X(\mathcal{C})$ is bounded over all $X$'s in within this radius. This is comparable to results in the factorized setting, for instance Chen and Wainwright [8, Theorem 1], where the initialization point is similarly assumed to be within a radius that is smaller than $\sigma_r(\widehat{X})$ of the true solution $\widehat{X}$.

**Approximate projections** The projection to $\mathcal{C}$, $P_\mathcal{C}$, can be obtained using a singular value decomposition (SVD), where only the top $r$ singular values and singular vectors of the matrix are retained to compute the best rank $r$ approximation. Nonetheless, it can be expensive to compute the SVD of a dense $n \times m$ matrix. We next propose an approximate projection operator for this space to avoid the cost of a singular value decomposition on an $n \times m$ matrix at each iteration of projected gradient descent.

To construct $P_X$, we first define some notation: for any rank-$r$ matrix $X$, let $T_X$ be the tangent space of low-rank matrices at $X$, given by[6]

$$T_X = \Big\{UA^\top + BV^\top \ : \ U \in \mathbb{R}^{n \times r}, V \in \mathbb{R}^{m \times r} \text{ are orthonormal bases for the column and row span of } X;$$
$$A \in \mathbb{R}^{m \times r}, B \in \mathbb{R}^{n \times r} \text{ are any matrices}\Big\}. \quad (23)$$

(This tangent space has frequently been studied in the context of nuclear norm minimization, see for instance [3].) We then define $P_X$ by first projecting to $T_X$, then projecting to the rank-$r$ constraint, that is,
$$P_X(Z) = P_\mathcal{C}\left(P_{T_X}(Z)\right). \quad (24)$$

---

[6]For $X \in \mathcal{C}$ which is of rank strictly lower than $r$, we can define $T_X$ by taking $U \in \mathbb{R}^{n \times r}, V \in \mathbb{R}^{m \times r}$ to be any orthonormal matrices which contain the column and row span of $X$; this choice is not unique, but formally we assume that we have fixed some choice of space $T_X$ for each $X \in \mathcal{C}$.

While this approximate projection will introduce some small error into the update steps, thus slowing convergence somewhat, it comes with a potentially large benefit: the SVD computations are always carried out on low rank matrices. Specifically, defining $U, V$ to be orthonormal bases for column and row spans of $X$ as before, for any $Z \in \mathbb{R}^d$ we can write

$$P_{T_X}(Z) = UU^\top Z + ZVV^\top - UU^\top ZVV^\top = \underbrace{\begin{pmatrix} U & (\mathbf{I}_n - UU^\top)ZV \end{pmatrix}}_{n \times 2r} \cdot \underbrace{\begin{pmatrix} Z^\top U & V \end{pmatrix}^\top}_{m \times 2r},$$

which means that calculating $P_\mathcal{C}(P_{T_X}(Z))$ can be substantially faster than the exact projection $P_\mathcal{C}(Z)$ when the rank bound $r$ is small while dimensions $n, m$ are large. Once this projection is computed, we now have new row and column span matrices $U, V$ ready to use for the next iteration's approximate projection step.

Our next result shows that this family of operators satisfies the conditions needed for our convergence results to be applied.

**Lemma 8.** *Let $\mathcal{C} = \{X \in \mathbb{R}^{n \times m} : \operatorname{rank}(X) \leq r\}$, and define the family of operators $P_X$ as in* (24). *Let $\rho = \frac{\sigma_r(\widehat{X})}{4}$. Then the inner product condition* (17)*, the norm compatibility condition* (18)*, and the local continuity condition* (20)*, are satisfied with $\gamma^\mathrm{c} = \gamma^\mathrm{d} = \frac{6}{\sigma_r(\widehat{X})}$ and $\phi = 3$, with respect to norms $\|\cdot\| = \|\cdot\|_\mathrm{nuc}$ and $\|\cdot\|^* = \|\cdot\|_\mathrm{sp}$.*

We see that up to a constant, this matches the results in Lemma 7 for the exact projection $P_\mathcal{C}$, and so we can expect roughly comparable convergence behavior with these approximate projections, while at the same time gaining computational efficiency by avoiding large singular value decompositions. We will compare this approximate projection method to exact projection and factored projection empirically in Section 6.

## 5.2 Sparsity

In many applications in high-dimensional statistics, the signal of interest is believed to be sparse or approximately sparse. Using an $\ell_1$ penalty or constraint serves as a convex relaxation to the sparsity constraint,

$$\arg\min_x \{\mathsf{g}(x) + \lambda \|x\|_1\} \text{ or } \arg\min_x \{\mathsf{g}(x) : \|x\|_1 \leq c\}$$

(i.e. the Lasso method [29], in the case of a linear regression problem). The convex $\ell_1$ norm penalty shrinks many coefficients to zero, but also leads to undesirable shrinkage bias on the large coefficients of $x$. Optimization with hard sparsity constraints (e.g. the iterative hard thresholding method [17]), while sometimes prone to local minima, are known to be successful in many settings and provide an alternative to convex relaxations (like the $\ell_1$ penalty) which induce shrinkage bias on large coefficients.

The shrinkage problem can also be alleviated by turning to nonconvex regularization functions, including the $\ell_q$ "norm" for $q < 1$, $\|x\|_q^q = \sum_i |x_i|^q$, whose convergence properties are studied by e.g. [18, 7], as well as the SCAD penalty [12], the MCP penalty [31], and the adaptive Lasso / reweighted $\ell_1$ method [4], which is related to a nonconvex "log-$\ell_1$" penalty of the form

$$\mathrm{logL1}_\nu(x) = \sum_i \nu \log(1 + |x_i|/\nu). \tag{25}$$

Smaller values of $\nu > 0$ correspond to greater nonconvexity, while setting $\nu = \infty$ recovers the $\ell_1$ norm.

Loh and Wainwright [20] study the convergence properties of a gradient descent algorithm for the penalized optimization problem $\arg\min_x \{\mathsf{g}(x) + \lambda \mathrm{Pen}(x)\}$, where the regularizer takes the form

$$\mathrm{Pen}(x) = \sum_i \mathsf{p}(|x_i|)$$

where $\mathsf{p}(t)$ is nondecreasing and concave over $t \geq 0$, but its concavity is bounded and it has finite derivative as $t \searrow 0$. Essentially, this means that $\mathrm{Pen}(x)$ behaves like a nonconvex version of the $\ell_1$ norm, shrinking small coefficients to zero but avoiding heavy shrinkage on large coefficients; the SCAD, MCP, and log-$\ell_1$ penalties are all examples. (The $\ell_q$ norm for $q < 1$ does not fit these assumptions, however, due to infinite derivative for coordinates $x_i \to 0$.) Proximal gradient descent with a nonconvex penalty such as SCAD is also studied by Lewis and Wright [19, Section 2.5] in the context of prox-regular functions. In this work, we consider the constrained version of this optimization problem, namely $\arg\min_x \{g(x) : \text{Pen}(x) \leq c\}$.

**Nonconvex regularizers** The general nonconvex sparsity-inducing penalties studied by Loh and Wainwright [20] are required to satisfy the following conditions (changing their notation slightly):

$$\text{Pen}(x) = \sum_i \mathsf{p}(|x_i|) \text{ where } \begin{cases} \mathsf{p}(0) = 0 \text{ and } \mathsf{p} \text{ is nondecreasing,} \\ t \mapsto \mathsf{p}(t)/t \text{ is nonincreasing (i.e. } \mathsf{p} \text{ is concave),} \\ t \mapsto \mathsf{p}(t) + \frac{\mu}{2}t^2 \text{ is convex,} \\ \mathsf{p} \text{ is differentiable on } t > 0, \text{ with } \lim_{t \searrow 0} \mathsf{p}'(t) = 1. \end{cases} \quad (26)$$

The following result calculates the local concavity coefficients for $\mathcal{C} = \{x : \text{Pen}(x) \leq c\}$.

**Lemma 9.** *Suppose that* $\text{Pen}(x) = \sum_i \mathsf{p}(|x_i|)$ *where* $\mathsf{p}$ *satisfies conditions* (26)*. Then*

$$\begin{cases} \gamma_x(\mathcal{C}) \leq \frac{\mu/2}{\mathsf{p}'(x_{\min})}, & \text{if } \text{Pen}(x) = c, \\ \gamma_x(\mathcal{C}) = 0, & \text{if } \text{Pen}(x) < c, \end{cases}$$

*with respect to the norm* $\|\cdot\| = \|\cdot\|_1$ *and its dual* $\|\cdot\|^* = \|\cdot\|_\infty$*, where for any* $x \in \mathbb{R}^d \backslash \{0\}$ *we define* $x_{\min}$ *to be the magnitude of its smallest nonzero entry.*

As an example, consider the log-$\ell_1$ penalty (25), so that our constraint set is

$$\mathcal{C} = \left\{ x \in \mathbb{R}^d : \sum_i \nu \log(1 + |x_i|/\nu) \leq c \right\}.$$

In this case we have $\mathsf{p}(t) = \nu \log(1 + t/\nu)$, which satisfies Loh and Wainwright [20]'s conditions (26) with $\mu = \frac{1}{\nu}$, and we can calculate $\mathsf{p}'(t) = \frac{1}{1+t/\nu}$. Therefore the local concavity coefficients for points $x$ on the boundary of $\mathcal{C}$ are bounded as $\gamma_x(\mathcal{C}) \leq \frac{1+x_{\min}/\nu}{2\nu}$. In particular, taking a maximum over all $x \in \mathcal{C}$, we obtain $\gamma(\mathcal{C}) \leq \frac{e^{c/\nu}}{2\nu}$. We can also check the norm compatibility condition:

**Lemma 10.** *If* $\mathcal{C} = \{x \in \mathbb{R}^d : \sum_i \mathsf{p}(|x_i|) \leq c\}$ *where* $c > 0$ *and* $\mathsf{p} : [0, \infty) \to [0, \infty)$ *satisfies the conditions* (26)*, then the norm compatibility condition* (13) *is satisfied with* $\phi = \frac{1}{\mathsf{p}'(\mathsf{p}^{-1}(c))}$*.*

With these results in place, we would therefore expect good convergence for projected gradient descent algorithms over the nonconvex sparsity constraint $\text{Pen}(x) \leq c$, as long as the objective function $\mathsf{g}$ and the initialization point satisfy the appropriate assumptions.

### 5.2.1 Constraints versus penalities

In Loh and Wainwright [20], the nonconvex optimization problem $\min_x \{\mathsf{g}(x) + \lambda \text{Pen}(x)\}$ is studied with no assumptions about the initial point $x_0$. Instead, they give an assumption that the concavity in Pen must be outweighed by the (restricted) convexity in $\mathsf{g}$. In our work on the constrained form of this problem, $\min\{\mathsf{g}(x) : \text{Pen}(x) \leq c\}$, we instead rely heavily on initialization conditions, namely (14), for projected gradient descent to succeed. While our result, Lemma 5, gives some justification for the necessity of the initialization conditions for general constraint sets $\mathcal{C}$, here we consider the specific setting of nonconvex sparsity penalties, and offer a more direct comparison of projected and proximal gradient descent (solving the constrained or penalized forms of the problem, respectively).

Lemma 5 suggests that the condition
$$2\gamma_x(\mathcal{C})\|\nabla \mathsf{g}(x)\|^* < \alpha, \tag{27}$$
where $\alpha$ is the restricted strong convexity parameter for the loss $\mathsf{g}$, is to some extent necessary for the success of projected gradient descent to be assured; otherwise projected gradient descent may have $x$ as a "bad" stationary point. How does this condition relate to the proximal gradient descent algorithm for the penalized form? Specifically, suppose that $x$ is a stationary point of proximal gradient descent for some step size $\eta > 0$, namely,
$$x = \arg\min_y \left\{ \frac{1}{2}\|y - (x - \eta\nabla\mathsf{g}(x))\|_2^2 + \eta\lambda \mathrm{Pen}(y) \right\}.$$

By first-order optimality conditions, we must therefore have
$$0 \in \partial\left(\frac{1}{2}\|y - (x - \eta\nabla\mathsf{g}(x))\|_2^2 + \eta\lambda\mathrm{Pen}(y)\right)\bigg|_{y=x},$$

where the subdifferential $\partial \mathrm{Pen}(y)$ is defined coordinatewise,
$$\partial \mathsf{p}(|y_i|) = \begin{cases} \mathsf{p}'(|y_i|) \cdot \mathrm{sign}(y_i), & y_i \neq 0, \\ [-1, 1], & y_i = 0. \end{cases}$$

In other words,
$$x - (x - \eta\nabla\mathsf{g}(x)) + \eta\lambda \partial \mathrm{Pen}(x) \ni 0,$$
and so for all $i$, $(\nabla\mathsf{g}(x))_i \in -\lambda\partial\mathsf{p}(|x_i|) \subset [-\lambda, \lambda]$. Therefore, applying the bound on the concavity coefficients given in Lemma 9, we have
$$2\gamma_x(\mathcal{C})\|\nabla\mathsf{g}(x)\|_\infty \leq \frac{\mu\lambda}{\mathsf{p}'(x_{\min})}$$
for any stationary point $x$ of the proximal gradient descent algorithm. The assumptions for convergence of proximal gradient descent in Loh and Wainwright [20] require that $\mu\lambda \leq$ (constant) $\cdot \alpha$. Therefore, up to some constant, we see that the condition (27) is *automatically* satisfied by any stationary point of proximal gradient descent, but for projected gradient descent this condition can instead fail and this allows for "bad" stationary points. Of course, if $\lambda$ is too large, then proximal gradient descent can also fail to find the global minimum; however, if $\lambda$ is chosen appropriately then no initialization condition is needed, while for the constrained form, an initialization condition is apparently necessary regardless of the constraint value $c$.

To some extent, this suggests that projected and proximal gradient descent may have fundamentally different behavior in the nonconvex setting, contradicting the notion that working with a constraint or a regularizer should lead to the same results (up to issues of tuning).

### 5.3 Spheres, orthogonal groups, and orthonormal matrices

We next consider a constraint set given by $\mathcal{C} = \{X \in \mathbb{R}^{n \times r} : X^\top X = \mathbf{I}_r\}$, the space of all orthonormal $n \times r$ matrices. We can also consider a related set, $\mathcal{C} = \{X \in \mathbb{R}^{n \times n} : X^2 = X, \mathrm{rank}(X) = r\}$, the set of all rank-$r$ projection matrices (with the orthogonal group as a special case when $r = n$). These examples have a different flavor than the low rank and sparse settings considered above; while the previous examples effectively bound the complexity of the signal (by finding latent sparse or low-rank structure), here we are instead enforcing special properties, namely orthogonality and/or unit norm. Optimization problems over these types of constraint sets arise, for instance, in PCA type problems where we would like to find the best low-rank representation of a data set.

First, we consider $n \times r$ orthonormal matrices:

**Lemma 11.** *Let $\mathcal{C} = \{X \in \mathbb{R}^{n \times r} : X^\top X = \mathbf{I}_r\}$, the space of orthogonal $n \times r$ matrices. Then $\mathcal{C}$ has local concavity coefficients $\gamma_X(\mathcal{C}) = \frac{1}{2}$ with respect to $\|\cdot\| = \|\cdot\|_{\mathrm{nuc}}$ and dual norm $\|\cdot\|^* = \|\cdot\|_{\mathrm{sp}}$. The norm compatibility condition* (13) *holds with $\phi = 1$.*

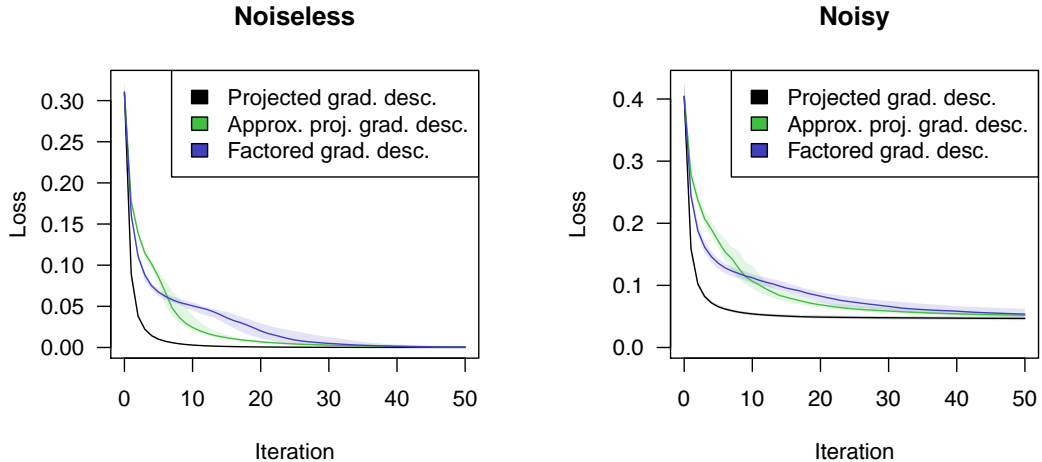

Figure 2: Results of the simulated matrix completion experiment comparing projected gradient descent on $\mathcal{C} = \{X \in \mathbb{R}^{100\times 100} : \mathrm{rank}(X) \leq 5\}$; gradient descent with approximate projections to the same set; and factored gradient descent on the variable $U \in \mathbb{R}^{100\times 5}$ (which relates to the low-rank matrix $X$ as $X = UU^\top$). For each method, the line and band show the median and quartiles of the loss $\mathrm{g}(X)$ over 50 trials.

Observe that the sphere $\mathbb{S}^{d-1} = \{x \in \mathbb{R}^d : \|x\|_2 = 1\}$ is a special case, obtained when $r = 1$.

Next, in many problems we may aim to find a rank-$r$ subspace that is optimal in some regard, but the exact choice of basis for this subspace does not matter; that is, an orthonormal basis $X \in \mathbb{R}^{n\times r}$ is identifiable only up to a rotation of its columns. In this case, we can instead choose to work with rank-$r$ projection matrices:

**Lemma 12.** *Let $\mathcal{C} = \{X \in \mathbb{R}^{n\times n} : \mathrm{rank}(X) = r, X \succeq 0, X^2 = X\}$, the space of rank-$r$ projection matrices. Then $\mathcal{C}$ has local concavity coefficients $\gamma_X(\mathcal{C}) \leq 2$ with respect to $\|\cdot\| = \|\cdot\|_{\mathrm{nuc}}$ and $\|\cdot\|^* = \|\cdot\|_{\mathrm{sp}}$.*

A special case is the setting $r = n$, when $\mathcal{C}$ is the orthogonal group.

# 6   Experiments

We next test the empirical performance of low rank optimization methods on a small matrix completion problem.[7] We generate an orthonormal matrix $U^\star \in \mathbb{R}^{100\times 5}$ at random, to produce a rank-5 positive semidefinite signal $U^\star U^{\star\top}$. We choose a subset of observed entries $\Omega \subset [100]\times [100]$ by giving each entry $(i,j)$ for $i \leq j$ a 20% chance of being observed, then symmetrizing across the diagonal. Our observations are then given by either $Y_{ij} = (U^\star U^{\star\top})_{ij}$ for the noiseless setting, or $Y_{ij} = (U^\star U^{\star\top})_{ij} + Z_{ij}$ for the noisy setting, where $Z_{ij} \sim N(0, \sigma^2)$ for entries $i \leq j$ and then we again symmetrize across the diagonal. We set $\sigma^2 = 0.2 \cdot \frac{\|U^\star U^{\star\top}\|_{\mathrm{F}}^2}{100^2}$ for a signal-to-noise ratio of 5. The loss function is

$$\mathrm{g}(X) = \frac{1}{2}\sum_{(i,j)\in\Omega}(Y_{ij} - X_{ij})^2,$$

and the constraint set is $\mathcal{C} = \{X \in \mathbb{R}^{100\times 100} : \mathrm{rank}(X) \leq 5\}$. We initialize at the matrix $X_0 = P_\mathcal{C}(Y_\Omega)$, where $Y_\Omega$ is the zero-filled matrix of observations, i.e. $(Y_\Omega)_{ij} = Y_{ij} \cdot \mathbb{1}\{(i,j)\in\Omega\}$.

We then compare three methods:

---

[7]Code for the simulated data experiments is available at http://www.stat.uchicago.edu/~rina/concavity.html.

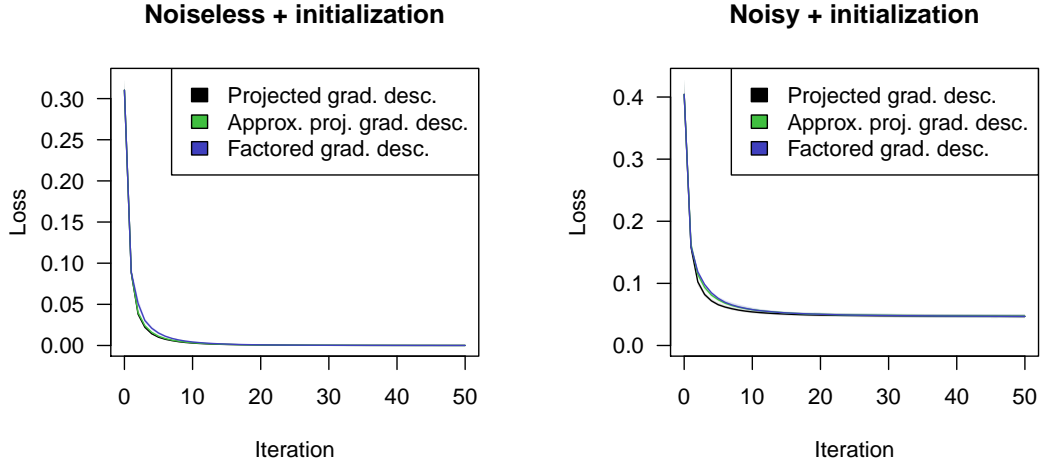

Figure 3: Same as in Figure 2, except that for all three methods, as an initialization step, the first iteration performs one step of exact projected gradient descent on $\mathcal{C}$; iterations $2, 3, \ldots, 50$ are then performed via the three different methods.

1. Projected gradient descent, where projection $P_\mathcal{C}$ to the rank constraint is carried out via a singular value decomposition at each step (i.e. the iterative hard thresholding method studied by Jain et al. [17]);

2. Approximate projected gradient descent, where the projection $P_\mathcal{C}$ is replaced by $P_X(Z) = P_\mathcal{C}(P_{T_X}(Z))$ as in (24), a more efficient computation;

3. Factored gradient descent (as studied by Chen and Wainwright [8], Zheng and Lafferty [33] for the positive semidefinite setting), where we define $\widetilde{\mathsf{g}}(U) = \mathsf{g}(UU^\top)$ over the variable $U \in \mathbb{R}^{100 \times 5}$, and perform (unconstrained) gradient descent on the variable $U$ with respect to $\widetilde{\mathsf{g}}$.

For each of the three methods, we run the method for 50 iterations with a constant step size $\eta$, and repeat for 50 trials. We then choose step size $\eta$ that achieves the lowest median loss at the last iteration, across all trials. Figure 2 displays the median loss, and the first and third quantiles, across the 50 trials, with respect to iteration number. We see that for both the noiseless and noisy settings, projected gradient descent achieves the lowest loss, with a very fast decay (while, of course, being the most computationally expensive method). Approximate projected gradient descent and factored gradient descent show an interesting comparison, where for early iterations ($\sim$5–10) the factored form gives a lower loss, while afterwards the approximate version performs better. By iteration 50, all three methods give nearly identical loss. It is likely that different settings may produce a different comparison between these methods.

Noting the rapid decrease in loss for the projected gradient descent method, we next ask whether the strengths of this method can be combined with the computational efficiency of the other two methods. As a second experiment, we repeat the steps above with the modification that, for all three methods being compared:

- Initialization step: for iteration 1, the update step is carried out with (exact) projected gradient descent;
- Then, for iterations $2, \ldots, 50$, the update step is carried out with the respective method.

The results are displayed in Figure 3. Remarkably, the three methods now perform nearly identically; starting with a single iteration of the more expensive projected gradient descent method is sufficient to allow the inexpensive methods to perform nearly as well.

# 7 Discussion

In this paper we have developed the local concavity coefficients, a measure of the extent to which a constraint set $\mathcal{C}$ violates convexity and may therefore be challenging for first-order optimization methods. These coefficients, related to the notion of prox-regularity in the analysis literature, are defined through four different measures of concavity that we then prove to be equivalent, connecting the geometric curvature of $\mathcal{C}$ with its behavior with respect to projections and first-order optimality conditions. This reveals a deep connection between geometry and optimization and allows us to analyze projected gradient descent to a range of examples such as low-rank estimation problems.

Many open questions remain in this area. As discussed earlier, the extent to which constrained versus penalized regularization (i.e. projected or proximal gradient descent) differ is not yet understood for nonconvex regularizers. In sparse estimation problems, the nonconvex $\ell_q$ "norm" is a popular regularizer that is empirically very successful (and has been studied theoretically), but it is not clear whether an $\ell_q$ norm constraint can fit into the framework of the concavity coefficients (i.e. whether $\gamma_x(\mathcal{C})$ is finite on the $\ell_q$ norm ball). For a low-rank estimation problem, research on factored gradient descent methods, which optimize over the function $U \mapsto \mathbf{g}(UU^\top)$ or $(U, V) \mapsto \mathbf{g}(UV^\top)$, has developed tools to work around the identifiability issue where factors are identifiable only up to rotation. Is there a more general way in which nonidentifiability, which can be thought of as a lack of convexity in certain directions, can be accounted for in the theory developed here?

Turning to our general results for convergence on an arbitrary nonconvex constraint set $\mathcal{C}$, it would be interesting to determine whether we can obtain a slower convergence rate assuming g is (restricted) Lipschitz and satisfies restricted strong convexity, but without a restricted smoothness result—standard results in the convex setting for this case suggest that we would want to take step size $\eta_t \propto 1/t$ and could expect a convergence rate of $\|x_t - \widehat{x}\|_2^2 \sim 1/t$. Finally, the strong initialization condition to ensure convergence to a global minimizer suggests that there may be some settings in which we can obtain weaker results—since our examples show that even convergence to a local minimum cannot be assured without checking the local concavity of the constraint set, are local concavity type assumptions sufficient to guarantee that projected gradient descent converges at least to some local minimizer?

# Acknowledgements

This work was partially supported by an Alfred P. Sloan fellowship and by NSF award DMS-1654076. The authors are grateful to John Lafferty and Stephen Wright for helpful feedback on an early version of this work.

## A  Proofs of local concavity coefficient results

In this section we prove the equivalence of the multiple notions of the (local or global) concavity of the constraint set $\mathcal{C}$, given in Theorems 1 and 2, as well as some properties of these coefficients (Lemmas 1, 2, 3, and 4).

Since our equivalent characterizations of the local concavity coefficients are inspired by many related conditions in the prox-regularity literature, some of the equivalence results are well-known in the $\ell_2$ setting (i.e. when the structured norm $\|\cdot\|$ is chosen to be the $\ell_2$ norm), as we have described in Section 2.3. Once we move to the general setting, where $\|\cdot\|$ may be any norm chosen to reflect the structure of the underlying signal, and where we do not assume that $x \mapsto \gamma_x(\mathcal{C})$ is continuous, many of the previously developed proof techniques will no longer apply. Throughout the proof, we will highlight those portions where our proof uses novel techniques due to the challenges of this more general setting.

In order to help discuss the various definitions of these coefficients before the equivalence is established, we begin by introducing notation for the local concavity coefficients defined using each of these four properties: for all $x \in \mathcal{C}$, define

$$\gamma_x^{\text{curv}}(\mathcal{C}) = \min\{\gamma \in [0, \infty] : \text{The curvature condition (1) holds for this point } x \text{ and any } y \in \mathcal{C}\},$$
$$\gamma_x^{\text{contr}}(\mathcal{C}) = \min\{\gamma \in [0, \infty] : \text{The contraction condition (3) holds for this point } x \text{ and any } y \in \mathcal{C}\},$$
$$\gamma_x^{\text{FO}}(\mathcal{C}) = \min\{\gamma \in [0, \infty] : \text{The first-order condition (4) holds for this point } x \text{ and any } y \in \mathcal{C}\},$$
$$\gamma_x^{\text{IP}}(\mathcal{C}) = \min\{\gamma \in [0, \infty] : \text{The inner product condition (5) holds for this point } x \text{ and any } y \in \mathcal{C}\}.$$

We emphasize that here we are *not* explicitly setting these coefficients to equal $\infty$ at degenerate points $x \in \mathcal{C}_{\text{dgn}}$— they may take finite values (we will need this distinction for some technical parts of our proofs later on). We will prove that these four definitions are all equal for all $x \notin \mathcal{C}_{\text{dgn}}$, which is sufficient for the equivalence result Theorem 2 since the local concavity coefficients are set to $\infty$ at degenerate points.

Before proceeding, we introduce one more definition: by equivalence of norms on $\mathbb{R}^d$, we can find some finite constant $B_{\text{norm}}$ such that

$$\text{For all } z \in \mathbb{R}^d, \begin{cases} B_{\text{norm}}^{-1}\|z\|_2 \leq \|z\| \leq B_{\text{norm}}\|z\|_2, \\ B_{\text{norm}}^{-1}\|z\|_2 \leq \|z\|^* \leq B_{\text{norm}}\|z\|_2. \end{cases} \quad (28)$$

Note that, while $B_{\text{norm}}$ is finite, it may be extremely large—for instance, $B_{\text{norm}} = \sqrt{d}$ when $\|\cdot\|$ is the $\ell_1$ norm.

## A.1 Proof of upper semi-continuity (Lemma 2)

Before we can prove the equivalence of the four definitions of the local coefficients in (6), we need to first show that these coefficients are upper semi-continuous, as claimed in Lemma 2. Of course, since we do not yet know that the four definitions are equivalent, we need to specify which definition we are using. We will work with the inner products property (5).

**Lemma 13.** *The map $x \mapsto \gamma_x^{\text{IP}}(\mathcal{C})$ is upper semi-continuous over $x \in \mathcal{C}\backslash\mathcal{C}_{\text{dgn}}$.*

This lemma will allow us to prove the equivalence result, Theorem 2. Once Theorem 2 is proved, then Lemma 13 becomes equivalent to the original lemma, Lemma 2, since $\gamma_x(\mathcal{C}) = \infty$ by definition on the subset $\mathcal{C}_{\text{dgn}} \subset \mathcal{C}$, which is a closed subset by definition, while Lemma 13 proves that $x \mapsto \gamma_x(\mathcal{C})$ is upper semi-continuous over the open subset $\mathcal{C}\backslash\mathcal{C}_{\text{dgn}} \subset \mathcal{C}$.

Before proving this result, we first state a well-known fact about projections, which we will use throughout our proofs:

$$\text{For any } z \in \mathbb{R}^d \text{ and } x \in \mathcal{C} \text{ with } P_\mathcal{C}(z) = x, \text{ for any } t \in [0, 1], P_\mathcal{C}((1-t)x + tz) = x. \quad (29)$$

Now we prove upper semi-continuity.

*Proof of Lemma 13.* Take any sequence $x_n \to x$, with $x, x_1, x_2, \dots \in \mathcal{C}\backslash\mathcal{C}_{\text{dgn}}$. We want to prove that

$$\gamma := \limsup_{n\to\infty} \gamma_{x_n}^{\text{IP}}(\mathcal{C}) \leq \gamma_x^{\text{IP}}(\mathcal{C}). \quad (30)$$

Since $x \notin \mathcal{C}_{\text{dgn}}$ by assumption, we know that $P_\mathcal{C}$ is continuous in some neighborhood of $x$. Let $r > 0$ be some radius so that $P_\mathcal{C}$ is continuous on $\mathbb{B}_*(x, r)$, where $\mathbb{B}_*(x, r)$ is the ball of radius $r$ around the point $x$ in the dual norm $\|\cdot\|^*$. Assume also that $\gamma > 0$, otherwise again the claim is trivial.

Taking a subsequence of the points $x_1, x_2, \dots$ if necessary, we can assume without loss of generality that

$$\gamma_{x_n}^{\text{IP}}(\mathcal{C}) \to \gamma.$$

Fix any $\epsilon > 0$ such that $\epsilon < \gamma$. For each $n$, by definition of the local concavity coefficient $\gamma_{x_n}^{\text{IP}}(\mathcal{C})$, there must exist some $y_n \in \mathcal{C}$ and some $z_n' \in \mathbb{R}^d$ with $P_\mathcal{C}(z_n') = x_n$, such that

$$\langle y_n - x_n, z_n' - x_n \rangle > \left(\gamma_{x_n}^{\text{IP}}(\mathcal{C}) - \epsilon\right) \|z_n' - x_n\|^* \|y_n - x_n\|_2^2. \quad (31)$$

Define

$$z_n = \begin{cases} z_n', & \text{if } \|z_n' - x_n\|^* \leq r/2, \\ x_n + (z_n' - x_n) \cdot \frac{r/2}{\|z_n' - x_n\|^*}, & \text{if } \|z_n' - x_n\|^* > r/2, \end{cases}$$

so that $\|z_n - x_n\|^* \leq r/2$. By (29), $P_\mathcal{C}(z_n) = x_n$. Furthermore, rescaling both sides of the inequality (31),

$$\langle y_n - x_n, z_n - x_n \rangle > \left(\gamma_{x_n}^{\text{IP}}(\mathcal{C}) - \epsilon\right) \|z_n - x_n\|^* \|y_n - x_n\|_2^2. \quad (32)$$

Since the left-hand side is bounded by $\|y_n - x_n\| \|z_n - x_n\|^*$, we see that

$$\|y_n - x_n\|_2^2 < \frac{\|y_n - x_n\|}{\gamma_{x_n}^{\text{IP}}(\mathcal{C}) - \epsilon} \leq \frac{\|y_n - x_n\|}{(\gamma - \epsilon)/2}$$

for all $n$ sufficiently large so that $\gamma_{x_n}^{\text{IP}}(\mathcal{C}) > \gamma - \frac{\gamma - \epsilon}{2}$. Therefore, since $\|y_n - x_n\| \leq B_{\text{norm}} \|y_n - x_n\|_2$ for some finite $B_{\text{norm}}$, then for all large $n$, $y_n$ lies in some ball of finite radius around $x$. The same is true for $z_n$ since $\|z_n - x_n\|^* \leq r/2$ by construction. Thus we can find a convergent subsequence, that is, $n_1, n_2, \ldots$ such that

$$\begin{cases} y_{n_i} \to y \text{ for some point } y, \\ z_{n_i} \to z \text{ for some point } z. \end{cases}$$

Since $\mathcal{C}$ is closed, we must have $y \in \mathcal{C}$. And, since $x_{n_i} \to x$, for sufficiently large $i$ we have $\|x_{n_i} - x\|^* \leq r/2$, so that $z_{n_i} \in \mathbb{B}_*(x, r)$. Since $P_\mathcal{C}$ is continuous on the ball $\mathbb{B}_*(x, r)$, then, $P_\mathcal{C}(z_{n_i}) = x_{n_i} \to x$ implies that we must have $P_\mathcal{C}(z) = x$. And,

$$\langle y - x, z - x \rangle = \lim_{i \to \infty} \langle y_{n_i} - x_{n_i}, z_{n_i} - x_{n_i} \rangle \geq \lim_{i \to \infty} \left( \gamma_{x_{n_i}}^{\text{IP}}(\mathcal{C}) - \epsilon \right) \|z_{n_i} - x_{n_i}\|^* \|y_{n_i} - x_{n_i}\|_2^2$$
$$= (\gamma - \epsilon) \|z - x\|^* \|y - x\|_2^2,$$

where the inequality applies (32) for each $n_i$. Therefore, $\gamma_x^{\text{IP}}(\mathcal{C}) \geq \gamma - \epsilon$. Since $\epsilon > 0$ was chosen to be arbitrarily small, this proves that $\gamma_x^{\text{IP}}(\mathcal{C}) \geq \gamma$, as desired. $\square$

## A.2  Proof of equivalence for local concavity (Theorem 2)

Now that we have established upper semi-continuity of $\gamma_x^{\text{IP}}(\mathcal{C})$ over $x \in \mathcal{C} \backslash \mathcal{C}_{\text{dgn}}$, we are ready to prove the equivalence of the local concavity coefficients.

Recall that if $x \in \mathcal{C}_{\text{dgn}}$ then $\gamma_x(\mathcal{C}) = \infty$ under all four definitions. Therefore, from this point on, we only need to show that

$$\gamma_x^{\text{curv}}(\mathcal{C}) = \gamma_x^{\text{contr}}(\mathcal{C}) = \gamma_x^{\text{IP}}(\mathcal{C}) = \gamma_x^{\text{FO}}(\mathcal{C}) \text{ for all } x \in \mathcal{C} \backslash \mathcal{C}_{\text{dgn}}.$$

In fact, we will also show that a weaker statement holds for all $x \in \mathcal{C}$ (i.e. without excluding degenerate points), namely

$$\gamma_x^{\text{IP}}(\mathcal{C}) \leq \min\{\gamma_x^{\text{curv}}(\mathcal{C}), \gamma_x^{\text{contr}}(\mathcal{C}), \gamma_x^{\text{FO}}(\mathcal{C})\} \text{ for all } x \in \mathcal{C}. \tag{33}$$

This additional bound will be useful later in our characterization of the degenerate points, when we prove Lemma 3.

### A.2.1  Inner products ⇒ First-order optimality

Fix any $u \in \mathcal{C} \backslash \mathcal{C}_{\text{dgn}}$. Let $f : \mathbb{R}^d \to \mathbb{R}$ be differentiable, and suppose that $u$ is a local minimizer of $f$ over $\mathcal{C}$. By Rockafellar and Wets [26, Theorem 6.12], this implies that $-\nabla f(u) \in N_\mathcal{C}(u)$, where $N_\mathcal{C}(u)$ is the normal cone to $\mathcal{C}$ at the point $u$ (see Rockafellar and Wets [26, Definition 6.3]). By Colombo and Thibault [10, (12)], we know that the normal cone can be obtained by a limit of proximal normal cones,

$$N_\mathcal{C}(u) = \lim_{x \in \mathcal{C}, x \to u} \sup \underbrace{\{w \in \mathbb{R}^d : P_\mathcal{C}(x + \epsilon \cdot w) = x \text{ for some } \epsilon > 0\}}_{\text{Proximal normal cone to } \mathcal{C} \text{ at } x}.$$

Therefore, we can find some sequences $u_1, u_2, \cdots \in \mathcal{C}$, $w_1, w_2, \cdots \in \mathbb{R}^d$, and $\epsilon_1, \epsilon_2, \cdots > 0$, such that $P_\mathcal{C}(u_n + \epsilon_n \cdot w_n) = u_n$ for all $n \geq 1$, with $u_n \to u$ and $w_n \to -\nabla f(u)$.

Now fix any $y \in \mathcal{C}$. By the inner product condition (5), for each $n \geq 1$,

$$\langle y - u_n, w_n \rangle = \langle y - u_n, (u_n + w_n) - u_n \rangle \leq \gamma_{u_n}^{\text{IP}}(\mathcal{C}) \|w_n\|^* \|y - u_n\|_2^2.$$

Taking limits on both sides, since $u_n \to u$ and $w_n \to -\nabla f(u)$,

$$\langle y - u, -\nabla f(u) \rangle \leq \left( \limsup_{t \to \infty} \gamma_{u_n}^{\text{IP}}(\mathcal{C}) \right) \cdot \|\nabla f(u)\|^* \|y - u\|_2^2.$$

Finally, recall that Lemma 13 proves that $x \mapsto \gamma_x^{\text{IP}}(\mathcal{C})$ is upper semi-continuous over $x \in \mathcal{C}\backslash\mathcal{C}_{\text{dgn}}$, and $\mathcal{C}_{\text{dgn}} \subset \mathcal{C}$ is a closed subset. Since $u \in \mathcal{C}\backslash\mathcal{C}_{\text{dgn}}$, we therefore have $u_n \in \mathcal{C}\backslash\mathcal{C}_{\text{dgn}}$ for all sufficiently large $t$, and therefore $\limsup_{t \to \infty} \gamma_{u_n}^{\text{IP}}(\mathcal{C}) \leq \gamma_u^{\text{IP}}(\mathcal{C})$. This proves that $\gamma_u^{\text{FO}}(\mathcal{C}) \leq \gamma_u^{\text{IP}}(\mathcal{C})$, as desired.

In fact, we can formulate a more general version of the first-order optimality condition:

For any Lipschitz continuous $\mathsf{f} : \mathbb{R}^d \to \mathbb{R}$ such that $x$ is a local minimizer of $\mathsf{f}$ over $\mathcal{C}$,
$$\langle y - x, v \rangle \geq -\gamma\|v\|^* \|y - x\|_2^2 \text{ for some } v \in \partial\mathsf{f}(x), \quad (34)$$

where $\partial\mathsf{f}(x)$ is the subdifferential to $\mathsf{f}$ at $x$ (see Rockafellar and Wets [26, Definition 8.3]). To see why (34) holds, Rockafellar and Wets [26, Theorem 8.15] guarantees that, since $\mathsf{f}$ is Lipschitz and $x$ is a local minimizer of $\mathsf{f}$ over the closed set $\mathcal{C}$, then we must have $-v \in N_\mathcal{C}(x)$ for some subgradient $v \in \partial\mathsf{f}(x)$.[8] The remainder of the proof is identical to the differentiable case treated above, with $v$ in place of $\nabla\mathsf{f}(x)$; this proves that, for any $x \in \mathcal{C}\backslash\mathcal{C}_{\text{dgn}}$ and any $y \in \mathcal{C}$, the stronger first-order optimality condition (34) holds with $\gamma = \gamma_x^{\text{IP}}(\mathcal{C})$.

Comparing to proofs of related conditions in the literature, we see that avoiding a continuity assumption on the map $x \mapsto \gamma_x$ means that the first-order optimality condition (4) does not follow immediately from the inner product condition (5); however, by first establishing upper semi-continuity of the map $x \mapsto \gamma_x^{\text{IP}}(\mathcal{C})$, the result follows.

### A.2.2   First-order optimality $\Rightarrow$ Inner products

This direction of the equivalence is immediate from the definitions of these two conditions. Fix any $x, y \in \mathcal{C}$ and $z \in \mathbb{R}^d$ with $P_\mathcal{C}(z) = x$. Consider the function $\mathsf{f}(w) = \frac{1}{2}\|w - z\|_2^2$, so that $P_\mathcal{C}(z) = x \in \mathcal{C}$ minimizes $\mathsf{f}$ over $\mathcal{C}$. We see trivially that $\nabla\mathsf{f}(x) = x - z$ and so

$$\langle y - x, z - x \rangle = \langle y - x, -\nabla\mathsf{f}(x)\rangle \leq \gamma_x^{\text{FO}}(\mathcal{C})\|\nabla\mathsf{f}(x)\|^*\|x - y\|_2^2 = \gamma_x^{\text{FO}}(\mathcal{C})\|z - x\|^*\|x - y\|_2^2.$$

Therefore $\gamma_x^{\text{IP}}(\mathcal{C}) \leq \gamma_x^{\text{FO}}(\mathcal{C})$ for all $x \in \mathcal{C}$, while previously we showed that the reverse inequality holds over $x \in \mathcal{C}\backslash\mathcal{C}_{\text{dgn}}$. Therefore, $\gamma_x^{\text{IP}}(\mathcal{C}) = \gamma_x^{\text{FO}}(\mathcal{C})$ for $x \in \mathcal{C}\backslash\mathcal{C}_{\text{dgn}}$.

### A.2.3   Curvature $\Rightarrow$ Inner products

Fix any $x, y \in \mathcal{C}$ and any $z \in \mathbb{R}^d$ with $P_\mathcal{C}(z) = x$. For all $t \in (0, 1)$, let $x_t = (1 - t)x + ty$, and choose

$$\widetilde{x}_t \in \arg\min_{x \in \mathcal{C}} \|x - x_t\| \text{ such that } \limsup_{t \searrow 0} \frac{\|\widetilde{x}_t - x_t\|}{t} \leq \gamma_x^{\text{curv}}(\mathcal{C})\|x - y\|_2^2,$$

as in the definition of $\gamma_x^{\text{curv}}(\mathcal{C})$. Fix any $\epsilon > 0$. Then for some $t_0 > 0$, for all $t < t_0$,

$$\frac{\|\widetilde{x}_t - x_t\|}{t} \leq \gamma_x^{\text{curv}}(\mathcal{C})\|x - y\|_2^2 + \epsilon.$$

Since $x = P_\mathcal{C}(z)$, this means that for all $t \in (0, 1)$,

$$\|z - x\|_2^2 \leq \|z - \widetilde{x}_t\|_2^2 = \|z - x_t\|_2^2 + \|\widetilde{x}_t - x_t\|_2^2 + 2\langle z - x_t, x_t - \widetilde{x}_t \rangle.$$

We can also calculate

$$\|z - x_t\|_2^2 = \|z - (1-t)x - ty\|_2^2 = \|z - x\|_2^2 - 2t\langle y - x, z - x\rangle + t^2\|x - y\|_2^2.$$

---

[8]More precisely, Rockafellar and Wets [26, Theorem 8.15] assumes only that $\mathsf{f}$ is proper and lower semi-continuous, but additionally requires the condition that $\partial^\infty \mathsf{f}(x) \cap \big(-N_\mathcal{C}(x)\big) = \{0\}$ (see Rockafellar and Wets [26, Chapter 8] for definitions). Since the horizon subdifferential $\partial^\infty \mathsf{f}(x)$ contains only the zero vector for any Lipschitz function $\mathsf{f}$, this condition must be satisfied once we assume that $\mathsf{f}$ is Lipschitz.

We rearrange terms to obtain

$$\langle y - x, z - x \rangle \leq \frac{1}{2t}\left(\|\widetilde{x}_t - x_t\|_2^2 + 2\langle z - x_t, x_t - \widetilde{x}_t\rangle + t^2\|x - y\|_2^2\right).$$

Recalling that $\|\cdot\|_2 \leq B_{\text{norm}}\|\cdot\|$ for some finite constant $B_{\text{norm}}$ by (28), we then have

$$\langle y - x, z - x\rangle \leq \frac{1}{2t}\left((B_{\text{norm}})^2\|\widetilde{x}_t - x_t\|^2 + 2\|z - x_t\|^*\|\widetilde{x}_t - \widetilde{x}\| + t^2\|x - y\|_2^2\right)$$

$$\leq \frac{1}{2t}\left((B_{\text{norm}})^2\left((\gamma_x^{\text{curv}}(\mathcal{C})\|x - y\|_2^2 + \epsilon)\cdot t\right)^2 + 2\|z - x_t\|^*(\gamma_x^{\text{curv}}(\mathcal{C})\|x - y\|_2^2 + \epsilon)\cdot t + t^2\|x - y\|_2^2\right)$$

$$= \|z - x_t\|^*(\gamma_x^{\text{curv}}(\mathcal{C})\|x - y\|_2^2 + \epsilon) + \frac{t}{2}\left((B_{\text{norm}})^2\left((\gamma_x^{\text{curv}}(\mathcal{C})\|x - y\|_2^2 + \epsilon)\right)^2 + \|x - y\|_2^2\right).$$

Taking a limit as $t$ approaches zero,

$$\langle y - x, z - x\rangle \leq (\gamma_x^{\text{curv}}(\mathcal{C})\|x - y\|_2^2 + \epsilon)\cdot\|z - x\|^*.$$

Since $\epsilon > 0$ was chosen to be arbitrarily small, therefore, $\gamma_x^{\text{IP}}(\mathcal{C}) \leq \gamma_x^{\text{curv}}(\mathcal{C})$, for any $x \in \mathcal{C}$.

To compare this portion of the proof to the existing literature, in the $\ell_2$ setting, this equivalence is proved in Colombo and Thibault [10, Theorem 14(q)]. Their proof relies on the fact that the projection operator $P_\mathcal{C}$ is taken with respect to the $\ell_2$ norm, and the curvature condition also seeks to bound the $\ell_2$ distance between the point $x_t$ and the set $\mathcal{C}$. In our more general setting, the curvature condition (1) works with the structured norm $\|\cdot\|$, while projections are still taken with respect to the $\ell_2$ norm, and so the same proof technique can no longer be applied (for instance, using this technique in a sparse problem with $\|\cdot\| = \|\cdot\|_1$, our results would suffer a factor of $B_{\text{norm}} = \sqrt{d}$ by converting from the $\ell_1$ norm to the $\ell_2$ norm). In our proof, a careful treatment of these various notions of distance allows for the bound to hold.

### A.2.4   Inner products $\Rightarrow$ Curvature

To prove the curvature condition, we will actually need to use the stronger form (34) of the first-order optimality condition—as proved in Appendix A.2.1, this condition holds with $\gamma = \gamma_x^{\text{IP}}(\mathcal{C})$ for all $x \in \mathcal{C}\backslash\mathcal{C}_{\text{dgn}}$.

Fix any $u \in \mathcal{C}\backslash\mathcal{C}_{\text{dgn}}$ and $y \in \mathcal{C}$. Let $u_t = (1-t)\cdot u + t\cdot y$, and define $\mathsf{f}(x) = \|x - u_t\|$. Note that $\mathsf{f}$ is a Lipschitz function. Since $\mathcal{C}$ is closed, and $\mathsf{f}$ is continuous and nonnegative, it must attain a minimum over $\mathcal{C}$, $x_t \in \arg\min_{x\in\mathcal{C}}\mathsf{f}(x)$. Since $\mathcal{C}_{\text{dgn}}$ is a closed subset of $\mathcal{C}$, this means that $x_t \in \mathcal{C}\backslash\mathcal{C}_{\text{dgn}}$ for any sufficiently small $t > 0$, since

$$\|x_t - u\| \leq \|x_t - u_t\| + \|u_t - u\| \leq 2\|u - u_t\| = 2t\|u - y\|$$

(where the second inequality uses the definition of $x_t$), and so $x_t \to u$.

Next, consider the subdifferential $\partial\mathsf{f}(x_t)$. It is well known that this subdifferential is not empty, and any element $v \in \partial\mathsf{f}(x_t)$ must satisfy $\|v\|^* \leq 1$ and $\langle v, x_t - u_t\rangle = \|x_t - u_t\|$. Now, applying the stronger form of the first-order optimality condition given in (34), we have

$$\langle v, y - x_t\rangle \geq -\gamma_{x_t}^{\text{IP}}(\mathcal{C})\|v\|^*\|y - x_t\|_2^2 = -\gamma_{x_t}^{\text{IP}}(\mathcal{C})\|y - x_t\|_2^2$$

and similarly, replacing $y \in \mathcal{C}$ with $u \in \mathcal{C}$,

$$\langle v, u - x_t\rangle \geq -\gamma_{x_t}^{\text{IP}}(\mathcal{C})\|u - x_t\|_2^2.$$

Taking the appropriate linear combination of these two inequalities,

$$\langle v, x_t - u_t\rangle \leq \gamma_{x_t}^{\text{IP}}(\mathcal{C})\left((1-t)\|u - x_t\|_2^2 + t\|y - x_t\|_2^2\right) = \gamma_{x_t}^{\text{IP}}(\mathcal{C})\left(t(1-t)\|u - y\|_2^2 + \|u_t - x_t\|_2^2\right),$$

where the last step simply uses the definition $u_t = (1-t)u + ty$ and rearranges terms. Finally, $\|u_t - x_t\|_2 \leq B_{\text{norm}}\|u_t - x_t\| \leq B_{\text{norm}}\|u - u_t\| = tB_{\text{norm}}\|u - y\|$, by definition of $u_t$ and $x_t$, so combining everything we can write

$$\min_{x\in\mathcal{C}}\|x - u_t\| = \|x_t - u_t\| = \langle v, x_t - u_t\rangle \leq \gamma_{x_t}^{\text{IP}}(\mathcal{C})\left(t(1-t)\|u - y\|_2^2 + t^2 B_{\text{norm}}^2\|u - y\|^2\right).$$

Dividing by $t$ and taking a limit,

$$\lim_{t \searrow 0} \frac{\min_{x \in \mathcal{C}} \|x - u_t\|}{t} \leq \left( \limsup_{t \searrow 0} \gamma^{\text{IP}}_{x_t}(\mathcal{C}) \right) \cdot \|u - y\|_2^2.$$

Finally, recall that $x \mapsto \gamma^{\text{IP}}_x(\mathcal{C})$ is upper semi-continuous by Lemma 13, and $x_t \to u$ as proved above. We thus have $\limsup_{t \searrow 0} \gamma^{\text{IP}}_{x_t}(\mathcal{C}) \leq \gamma^{\text{IP}}_u(\mathcal{C})$. This proves that $\gamma^{\text{curv}}_u(\mathcal{C}) \leq \gamma^{\text{IP}}_u(\mathcal{C})$, for any $u \in \mathcal{C} \backslash \mathcal{C}_{\text{dgn}}$.

Combining with our previous steps, we now have

$$\gamma^{\text{FO}}_x(\mathcal{C}) = \gamma^{\text{IP}}_x(\mathcal{C}) = \gamma^{\text{curv}}_x(\mathcal{C})$$

for all $x \in \mathcal{C} \backslash \mathcal{C}_{\text{dgn}}$, while for $x \in \mathcal{C}_{\text{dgn}}$ we have the weaker statement $\gamma^{\text{IP}}_x(\mathcal{C}) \leq \min\{\gamma^{\text{curv}}_x(\mathcal{C}), \gamma^{\text{FO}}_x(\mathcal{C})\}$.

To compare this portion of the proof to the existing literature, in the $\ell_2$ setting, the equivalent result is proved in Colombo and Thibault [10, Proposition 9]. In their proof, they use the identity $\|x_t - u_t\|_2 = \frac{\langle x_t - u_t, x_t - u_t \rangle}{\|x_t - u_t\|_2}$, and then upper bound the right-hand side via the inner product condition. To translate this proof into the more general structured norm setting, we write $\|x_t - u_t\| = \langle v, x_t - u_t \rangle$ for a subgradient $v$ in the subdifferential of the function $x \mapsto \|x - u_t\|$, and apply results from the analysis literature along with our upper semi-continuity result, Lemma 13.

### A.2.5   Approximate contraction $\Leftrightarrow$ Inner products

This proof, for the case of a general norm $\|\cdot\|$, proceeds identically as the proof for the case where $\|\cdot\| = \|\cdot\|_2$ (presented e.g. in Colombo and Thibault [10, Theorem 3(b,d)]). For completeness, we reproduce the argument here.

First, we show that $\gamma^{\text{IP}}_x(\mathcal{C}) \leq \gamma^{\text{contr}}_x(\mathcal{C})$. Fix any $x \in \mathcal{C}$, and any $z \in \mathbb{R}^d$ with $x = P_\mathcal{C}(z)$. Define $z_t = t \cdot z + (1-t) \cdot x$ for $t \in [0, 1]$. By (29), $x = P_\mathcal{C}(z_t)$ for all $t \in [0, 1]$.

Then for any $y \in \mathcal{C}$, since $\|z_t - x\|^* = t\|z - x\|^*$,

$$\|y - x\|_2 \left( 1 - \gamma^{\text{contr}}_x(\mathcal{C}) \cdot t \|z - x\|^* \right) \leq \|y - z_t\|_2$$

by the approximate contraction property (3). For sufficiently small $t$, the left-hand side is nonnegative (except for the trivial case $\gamma^{\text{contr}}_x(\mathcal{C}) = \infty$, in which case there is nothing to prove). Squaring both sides and rearranging some terms,

$$\|y - x\|_2^2 \leq \|y - z_t\|_2^2 + (2\gamma^{\text{contr}}_x(\mathcal{C}) \cdot t\|z - x\|^* - (\gamma^{\text{contr}}_x(\mathcal{C}) \cdot t\|z - x\|^*)^2)\|y - x\|_2^2.$$

And,

$$\|y - z_t\|_2^2 = \|y - x\|_2^2 + \|x - z_t\|_2^2 + 2\langle y - x, x - z_t \rangle$$

so rearranging terms again,

$$2\langle y - x, z_t - x \rangle \leq \|x - z_t\|_2^2 + (2\gamma^{\text{contr}}_x(\mathcal{C}) \cdot t\|z - x\|^* - (\gamma^{\text{contr}}_x(\mathcal{C}) \cdot t\|z - x\|^*)^2)\|y - x\|_2^2.$$

Plugging in the definition of $z_t$,

$$2t\langle y - x, z - x \rangle \leq t^2\|x - z\|_2^2 + (2\gamma^{\text{contr}}_x(\mathcal{C}) \cdot t\|z - x\|^* - \gamma^{\text{contr}}_x(\mathcal{C})^2 \cdot t^2(\|z - x\|^*)^2)\|y - x\|_2^2.$$

Dividing by $2t$, then taking the limit as $t \searrow 0$,

$$\langle y - x, z - x \rangle \leq \gamma^{\text{contr}}_x(\mathcal{C}) \|z - x\|^* \|y - x\|_2^2.$$

Therefore, for any $x \in \mathcal{C}$, $\gamma^{\text{IP}}_x(\mathcal{C}) \leq \gamma^{\text{contr}}_x(\mathcal{C})$.

Now we prove the reverse inequality, i.e. $\gamma^{\text{contr}}_x(\mathcal{C}) \leq \gamma^{\text{IP}}_x(\mathcal{C})$. Fix any $x, y \in \mathcal{C}$ and any $z \in \mathbb{R}^d$ with $x = P_\mathcal{C}(z)$. Then

$$\|y - x\|_2^2 + \langle y - x, z - y \rangle = \langle y - x, z - x \rangle \leq \gamma^{\text{IP}}_x(\mathcal{C})\|z - x\|^* \|y - x\|_2^2.$$

Simplifying,
$$\left(1 - \gamma_x^{\text{IP}}(\mathcal{C})\|z - x\|^*\right)\|y - x\|_2^2 \leq -\langle y - x, z - y\rangle \leq \|y - x\|_2\|z - y\|_2,$$
and so
$$\left(1 - \gamma_x^{\text{IP}}(\mathcal{C})\|z - x\|^*\right)\|y - x\|_2 \leq \|z - y\|_2.$$
Therefore, for any $x \in \mathcal{C}$ $\gamma_x^{\text{contr}}(\mathcal{C}) \leq \gamma_x^{\text{IP}}(\mathcal{C})$.

Combining everything, we have now proved
$$\gamma_x^{\text{contr}}(\mathcal{C}) = \gamma_x^{\text{IP}}(\mathcal{C}) = \gamma_x^{\text{FO}}(\mathcal{C}) = \gamma_x^{\text{curv}}(\mathcal{C})$$
for all $x \in \mathcal{C}\backslash\mathcal{C}_{\text{dgn}}$, in addition to the weaker bound (33) for all $x \in \mathcal{C}$, as desired. This completes the proof of Theorem 2.

## A.3    Proof of characterization of degenerate points (Lemma 3)

Next we prove that the degenerate points $u \in \mathcal{C}_{\text{dgn}}$ are precisely those points where any of the four local concavity conditions would fail to hold, in any neighborhood of $u$ and for any finite $\gamma$. First, the characterization of prox-regularity given in Poliquin et al. [25, Proposition 1.2, Theorem 1.3(i)] proves that, if the projection operator $P_\mathcal{C}$ is *not* continuous in a neighborhood of $u \in \mathcal{C}$, then there are no constants $\epsilon > 0$ and $\gamma < \infty$ such that the inner product condition (5) holds for all $x \in \mathcal{C} \cap \mathbb{B}_2(u, \epsilon)$. Therefore, for any $r > 0$, $\sup_{x \in \mathcal{C} \cap \mathbb{B}_2(u,r)} \gamma_x^{\text{IP}}(\mathcal{C}) = \infty$.

Finally, in proving Theorem 2, we proved (33), i.e. $\gamma_x^{\text{IP}}(\mathcal{C}) \leq \min\{\gamma_x^{\text{curv}}(\mathcal{C}), \gamma_x^{\text{contr}}(\mathcal{C}), \gamma_x^{\text{FO}}(\mathcal{C})\}$ for all $x \in \mathcal{C}$. This implies that,
$$\lim_{r \to 0}\left\{\sup_{x \in \mathcal{C} \cap \mathbb{B}_2(u,r)} \gamma_x^{(*)}(\mathcal{C})\right\} = \infty,$$
where (*) denotes any of the four properties, i.e. $\gamma_x^{\text{curv}}(\mathcal{C})$ for the curvature condition (1), $\gamma_x^{\text{contr}}(\mathcal{C})$ for the contraction property (3), $\gamma_x^{\text{IP}}(\mathcal{C})$ for the inner product condition (5), or $\gamma_x^{\text{FO}}(\mathcal{C})$ for the first-order optimality condition (4). This proves the lemma.

## A.4    Proof of two-sided contraction property (Lemma 4)

This proof, for the case of a general norm $\|\cdot\|$, proceeds identically as the proof for the case where $\|\cdot\| = \|\cdot\|_2$ (presented e.g. in Colombo and Thibault [10, Theorem 3(b,d)]). For completeness, we reproduce the argument here.

Take any $x, y \in \mathcal{C}$ and any $z, w \in \mathbb{R}^d$ with $P_\mathcal{C}(z) = x$ and $P_\mathcal{C}(w) = x$. By definition of the local concavity coefficients, applying the inner product bound (5) we have
$$\langle y - x, z - x\rangle \leq \gamma_x(\mathcal{C})\|z - x\|^*\|x - y\|_2^2.$$
Applying the same property with the roles of the variables reversed,
$$\langle x - y, w - y\rangle \leq \gamma_y(\mathcal{C})\|w - y\|^*\|x - y\|_2^2.$$
Adding these two inequalities together,
$$\langle y - x, z - x - w + y\rangle \leq \gamma_x(\mathcal{C})\|z - x\|^*\|x - y\|_2^2 + \gamma_y(\mathcal{C})\|w - y\|^*\|x - y\|_2^2.$$
Rearranging terms and simplifying,
$$\left(1 - \gamma_x(\mathcal{C})\|z - x\|^* - \gamma_y(\mathcal{C})\|w - y\|^*\right)\|x - y\|_2^2 \leq \langle y - x, w - z\rangle.$$
Since the right-hand side is bounded by $\|x - y\|_2\|z - w\|_2$ by the Cauchy–Schwarz inequality, this proves the lemma.

In fact, we will also prove a related result that will be useful for our convergence proofs later on. As above, we have
$$\langle y - x, z - x \rangle \leq \gamma_x(\mathcal{C})\|z - x\|^* \|x - y\|_2^2,$$
and since $y = P_\mathcal{C}(w)$, we also have
$$0 \geq \|y - w\|_2^2 - \|x - w\|_2^2 = \|x - y\|_2^2 + 2\langle y - x, x - w \rangle.$$
Then, adding the two bounds together,
$$\langle y - x, z - w \rangle = \langle y - x, z - x \rangle + \langle y - x, x - w \rangle \leq \gamma_x(\mathcal{C})\|z - x\|^* \|x - y\|_2^2 - \frac{1}{2}\|x - y\|_2^2,$$
and so
$$\left(\frac{1}{2} - \gamma_x(\mathcal{C})\|z - x\|^*\right)\|x - y\|_2^2 \leq -\langle y - x, z - w \rangle \leq \|x - y\|_2 \|z - w\|_2.$$
This proves that
$$\left(\frac{1}{2} - \gamma_x(\mathcal{C})\|z - x\|^*\right)\|x - y\|_2 \leq \|z - w\|_2. \tag{35}$$
In a setting where $\gamma_x(\mathcal{C})\|z - x\|^*$ is small but $\gamma_y(\mathcal{C})\|w - y\|^*$ may be large, this alternate result can give a stronger bound than Lemma 4.

## A.5   Proof of equivalence for global concavity (Theorem 1) and local vs global coefficients (Lemma 1)

We prove Theorem 1, which states that the five definitions for the global concavity coefficient $\gamma(\mathcal{C})$ are equivalent, alongside Lemma 1, which states that $\gamma(\mathcal{C}) = \sup_{x \in \mathcal{C}} \gamma_x(\mathcal{C})$.

First, suppose that $\mathcal{C}$ contains one or more degenerate points, $\mathcal{C}_{\mathsf{dgn}} \neq \varnothing$, in which case $\sup_{x \in \mathcal{C}} \gamma_x(\mathcal{C}) = \infty$. By definition of $\mathcal{C}_{\mathsf{dgn}}$, the projection operator $P_\mathcal{C}$ is not continuous on any neighborhood of $\mathcal{C}$. Poliquin et al. [25, Theorem 4.1] prove that this implies $\mathcal{C}$ is not prox-regular, and so $\gamma(\mathcal{C}) = \infty$ as discussed in Section 2.3.

Next, suppose that $\mathcal{C}$ contains no degenerate points. Let $\gamma^* = \sup_{x \in \mathcal{C}} \gamma_x(\mathcal{C})$. Then clearly, by definition of the local coefficients $\gamma_x(\mathcal{C})$,
$$\gamma^* = \min\{\gamma \in [0, \infty] : \text{Property (*) holds for all } x, y \in \mathcal{C}\}$$
where (*) may refer to any of the four equivalent properties, namely the curvature condition (1), the (one-sided) contraction property (3), the inner product condition (5), and the first-order condition (4). Next, let
$$\gamma^\sharp = \min\{\gamma \in [0, \infty] : \text{The two-sided contraction property (2) holds for all } x, y \in \mathcal{C}\}.$$

Clearly, the two-sided contraction property (2) is stronger than its one-sided version (3), and so $\gamma^* \leq \gamma^\sharp$. However, Lemma 4 shows that they are in fact equal, proving that
$$\left(1 - \gamma_x(\mathcal{C})\|z - x\|^* - \gamma_y(\mathcal{C})\|w - y\|^*\right) \cdot \|x - y\|_2 \leq \|z - w\|_2$$
for all $z, w \in \mathbb{R}^d$ with $x = P_\mathcal{C}(z), y = P_\mathcal{C}(w)$. Since $\gamma_x(\mathcal{C}), \gamma_y(\mathcal{C}) \leq \gamma^*$ for all $x, y \in \mathcal{C}$, this implies that
$$\left(1 - \gamma^*\|z - x\|^* - \gamma^*\|w - y\|^*\right) \cdot \|x - y\|_2 \leq \|z - w\|_2,$$
that is, (2) holds for all $x, y \in \mathcal{C}$ with constant $\gamma = \gamma^*$. So, we have $\gamma^\sharp \leq \gamma^*$. Therefore, the five conditions defining $\gamma(\mathcal{C})$ are equivalent, and $\gamma(\mathcal{C}) = \gamma^\sharp = \gamma^* = \sup_{x \in \mathcal{C}} \gamma_x(\mathcal{C})$, proving Theorem 1 and Lemma 1.

# B   Proofs of convergence results

In this section we prove our convergence results for projected gradient descent (Theorem 3) and approximate projected gradient descent (Theorem 4), along with the necessity of the initialization condition (Lemma 5) and equivalence of the exact and approximate convergence results (Lemma 6).

## B.1   Proof of Theorem 3

This result, using the exact projection operator $P_\mathcal{C}$, is in fact a special case of Theorem 4, which provides a convergence guarantee using a family of operators $\{P_x\}$. To see why, define $P_x = P_\mathcal{C}$ for all $x \in \mathcal{C}$. To apply Theorem 4, we only need to check that the relevant assumptions, namely (17), (18), (19), and (20), all hold.

To check (17), by setting
$$\gamma^c = \max\{\gamma_x(\mathcal{C}) : x \in \mathcal{C} \cap \mathbb{B}_2(\widehat{x}, \rho)\} \text{ and } \gamma^d = 0,$$
we see that the desired bound is a trivial consequence of the inner product condition (5). The norm compatibility condition (18) for $\{P_x\}$ holds as a trivial consequence of the original norm compatibility condition (13). The initialization condition (19) follows directly from the original initialization condition (14) by our choice of $\gamma^c, \gamma^d$. Finally, we verify the local continuity condition (20). Fix any $x \in \mathcal{C} \cap \mathbb{B}_2(\widehat{x}, \rho)$ and any $z \in \mathbb{R}^d$ such that $P_\mathcal{C}(z) \in \mathbb{B}_2(\widehat{x}, \rho)$ and $2(\gamma^c + \gamma^d)\|z - P_x(z)\|^* \leq 1 - c_0$. By (35), for all $w \in \mathbb{R}^d$,
$$\left(\frac{1}{2} - \gamma_{P_\mathcal{C}(z)}(\mathcal{C})\|z - P_\mathcal{C}(z)\|^*\right) \|P_\mathcal{C}(z) - P_\mathcal{C}(w)\|_2 \leq \|z - w\|_2.$$

Since $P_x(z) = P_\mathcal{C}(z) \in \mathbb{B}_2(\widehat{x}, \rho)$ and $\gamma^c + \gamma^d = \max_{u \in \mathcal{C} \cap \mathbb{B}_2(\widehat{x}, \rho)} \gamma_u(\mathcal{C}) \geq \gamma_{P_\mathcal{C}(z)}(\mathcal{C})$, then,
$$\|P_\mathcal{C}(z) - P_\mathcal{C}(w)\|_2 \leq \frac{\|z - w\|_2}{1/2 - (\gamma^c + \gamma^d)\|z - P_x(z)\|^*} \leq \frac{\|z - w\|_2}{c_0/2}.$$

Setting $\delta = \epsilon \cdot c_0/2$ then proves the condition (20).

With these conditions in place, Theorem 3 becomes simply a special case of Theorem 4.

## B.2   Proof of Theorem 4

For $t = 0$, the statement holds trivially. To prove that the bound holds for subsequent steps, we will proceed by induction. Choose any $\rho_0 \in (0, \rho)$ such that
$$\rho_0 \geq \max\left\{\|x_0 - \widehat{x}\|_2, \sqrt{\frac{1.5\varepsilon_g^2}{c_0}}\right\},$$
where this maximum is $< \rho$ by assumption of the theorem. We will prove that
$$\begin{cases} \|x_{t+1} - \widehat{x}\|_2^2 \leq \left(1 - \frac{2c_0\alpha}{\text{Denom}}\right) \|x_t - \widehat{x}\|_2^2 + \frac{3\alpha}{\text{Denom}}\varepsilon_g^2, \\ \|x_{t+1} - \widehat{x}\|_2 \leq \rho_0, \end{cases} \tag{36}$$
for all $t \geq 0$, where the denominator term is given by
$$\text{Denom} = \beta + \alpha\left(c_0 - (1 - c_0) \cdot \frac{\gamma^c}{\gamma^c + \gamma^d}\right) \leq \alpha + \beta.$$

Assuming that this holds, then applying the first bound of (36) iteratively, we will then have
$$\|x_t - \widehat{x}\|_2^2 \leq \left(1 - \frac{2c_0\alpha}{\text{Denom}}\right)^t \|x_0 - \widehat{x}\|_2^2 + \frac{1.5}{c_0}\varepsilon_g^2 \leq \left(1 - \frac{2c_0\alpha}{\alpha + \beta}\right)^t \|x_0 - \widehat{x}\|_2^2 + \frac{1.5}{c_0}\varepsilon_g^2,$$

which proves the theorem.

Now we turn to proving (36), assuming that it holds at the previous time step. In order to apply assumptions such as restricted strong convexity (11) and the inner product condition (17), we first need to know that $\|x_{t+1} - \widehat{x}\|_2 \leq \rho$, which we cannot ensure directly at the start. To get around this, we first consider reducing the step size. For any step size $s \in [0, \eta]$, define

$$x'_{t+1}(s) = x_t - s\nabla \mathsf{g}(x_t) \text{ and } x_{t+1}(s) = P_{x_t}(x'_{t+1}(s)).$$

Define

$$\mathcal{S} = \{s \in [0, \eta] : \|x_{t+1}(s) - \widehat{x}\|_2 \leq \rho\} \text{ and } \mathcal{S}_0 = \{s \in [0, \eta] : \|x_{t+1}(s) - \widehat{x}\| \leq \rho_0\}.$$

Clearly $0 \in \mathcal{S}_0 \subseteq \mathcal{S}$ since $x_{t+1}(0) = x_t$, which satisfies $\|x_t - \widehat{x}\|_2 \leq \rho_0$ by assumption. First, we claim that we can find some $\Delta > 0$ such that

$$\text{If } s \in \mathcal{S}_0 \text{ then } \min\{s + \Delta, \eta\} \in \mathcal{S}. \tag{37}$$

To prove this, we will apply the local uniform continuity assumption (20). Let $x = x_t$ and $\epsilon = \rho - \rho_0$, and find $\delta > 0$ as in the assumption (20). For any $s \in \mathcal{S}_0$, let $z = x'_{t+1}(s) = x_t - s\nabla\mathsf{g}(x_t)$. Then $P_x(z) = P_{x_t}(x'_{t+1}(s)) = x_{t+1}(s) \in \mathbb{B}_2(\widehat{x}, \rho_0) \subset \mathbb{B}_2(\widehat{x}, \rho)$, and

$$2(\gamma^c + \gamma^d)\|z - P_x(z)\|^* \leq 2(\gamma^c + \gamma^d) \cdot \phi\|z - x_t\|^* = 2(\gamma^c + \gamma^d)\phi s\|\nabla\mathsf{g}(x_t)\|^*$$
$$\leq 2(\gamma^c + \gamma^d) \cdot \phi\eta \max_{x \in \mathcal{C} \cap \mathbb{B}_2(\widehat{x}, \rho)} \|\nabla\mathsf{g}(x)\|^* \leq \eta \cdot (1-c_0)\alpha \leq 1 - c_0,$$

where the first inequality uses the norm compatibility condition (18), the third uses the initialization condition (19), and the fourth uses $\eta = 1/\beta \leq 1/\alpha$. Therefore, the conditions of the local continuity statement (20) are satisfied, and so for any $w \in \mathbb{R}^d$ with $\|w - z\|_2 \leq \delta$, we must have $\|P_x(w) - P_x(z)\|_2 \leq \epsilon = \rho - \rho_0$. Then

$$\|P_x(w) - \widehat{x}\|_2 \leq \|P_x(w) - P_x(z)\|_2 + \|P_x(z) - \widehat{x}\|_2 \leq (\rho - \rho_0) + \rho_0 = \rho.$$

Now, define $\Delta = \delta/\|\mathsf{g}(x_t)\|_2$ and set $w = x'_{t+1}(\min\{s + \Delta, \eta\})$. Then $\|w - x'_{t+1}(s)\|_2 \leq \Delta\|\nabla\mathsf{g}(x_t)\|_2 \leq \delta$, and so $\|P_x(w) - \widehat{x}\|_2 \leq \rho$. This proves that $\min\{s + \Delta, \eta\} \in \mathcal{S}$, and so (37) holds.

Next, consider any $s \in \mathcal{S}$. If $s = 0$, then $s \in \mathcal{S}_0$ since $\|x_t - \widehat{x}\|_2 \leq \rho_0$ by the previous step. Otherwise, assume $s > 0$. We first have

$$\max\{\|x'_{t+1}(s) - x_{t+1}(s)\|^*, \|x'_{t+1}(s) - x_t\|^*\} = \max\{\|x'_{t+1}(s) - P_{x_t}(x'_{t+1}(s))\|^*, \|x'_{t+1}(s) - x_t\|^*\}$$
$$\leq \phi\|x'_{t+1}(s) - x_t\|^* = \phi\|-s\nabla\mathsf{g}(x_t)\|^* \leq \frac{s\alpha(1-c_0)}{2(\gamma^c + \gamma^d)}, \tag{38}$$

where first inequality uses the norm compatibility condition (18) while the second uses the initialization condition (19), since $\|x_t - \widehat{x}\|_2 \leq \rho$.

We can now use the inner product condition (17), applied with $x = x_t$ and $z = x'_{t+1}(s)$, with $P_x(z) = P_{x_t}(x'_{t+1}(s)) = x_{t+1}(s)$. (This condition can be applied as we have checked that $x_t, x_{t+1}(s) \in \mathbb{B}_2(\widehat{x}, \rho)$.) The inner product condition yields

$$\langle \widehat{x} - x_{t+1}(s), x'_{t+1}(s) - x_{t+1}(s) \rangle$$
$$\leq \max\{\|x'_{t+1}(s) - x_{t+1}(s)\|^*, \|x'_{t+1}(s) - x_t\|^*\} \cdot \left(\gamma^c\|\widehat{x} - x_{t+1}(s)\|_2^2 + \gamma^d\|\widehat{x} - x_t\|_2^2\right)$$
$$\leq \frac{s\alpha(1-c_0)}{2(\gamma^c + \gamma^d)}\left(\gamma^c\|\widehat{x} - x_{t+1}(s)\|_2^2 + \gamma^d\|\widehat{x} - x_t\|_2^2\right), \tag{39}$$

where the last step applies (38).

Next, abusing notation, define for every $u \in \mathcal{C}$

$$\gamma_u^c(\mathcal{C}) = \begin{cases} \gamma^c, & \|u - \widehat{x}\|_2 < \rho, \\ \infty, & \|u - \widehat{x}\|_2 \geq \rho, \end{cases} \text{ and } \gamma_u^d(\mathcal{C}) = \begin{cases} \gamma^d, & \|u - \widehat{x}\|_2 < \rho, \\ \infty, & \|u - \widehat{x}\|_2 \geq \rho. \end{cases}$$

Trivially these maps are upper semi-continuous. We also need to check the local continuity assumption (22), which requires that if $\gamma_x^c(\mathcal{C}) + \gamma_x^d(\mathcal{C}) < \infty$ and $z_t \to x$, then $P_x(z_t) \to x$. In fact, by the norm compatibility assumption (18), for all $x \in \mathcal{C} \cap \mathbb{B}_2(\widehat{x}, \rho)$ we have $\|z_t - P_x(z_t)\|^* \leq \phi \|z_t - x\|^*$, and so

$$\|P_x(z_t) - x\|^* \leq \|P_x(z_t) - z_t\|^* + \|z_t - x\|^* \leq (1 + \phi)\|z_t - x\|^* \to 0,$$

proving that $P_x(z_t) \to x$, as desired. By Lemma 6, the local concavity coefficients of $\mathcal{C}$ then satisfy $\gamma_u(\mathcal{C}) \leq \gamma_u^c(\mathcal{C}) + \gamma_u^d(\mathcal{C})$, and so in particular, $\gamma_{\widehat{x}}(\mathcal{C}) \leq \gamma^c + \gamma^d$.

We will now apply the first-order optimality conditions (4) at the point $x = \widehat{x}$. We have

$$\mathsf{g}(x_{t+1}(s)) - \mathsf{g}(\widehat{x}) \geq \langle x_{t+1}(s) - \widehat{x}, \nabla \mathsf{g}(\widehat{x}) \rangle + \frac{\alpha}{2}\|x_{t+1}(s) - \widehat{x}\|_2^2 - \frac{\alpha \varepsilon_\mathsf{g}^2}{2} \text{ by restricted strong convexity (11)}$$

$$\geq -(\gamma^c + \gamma^d)\|\nabla \mathsf{g}(\widehat{x})\|^* \|x_{t+1}(s) - \widehat{x}\|_2^2 + \frac{\alpha}{2}\|x_{t+1}(s) - \widehat{x}\|_2^2 - \frac{\alpha \varepsilon_\mathsf{g}^2}{2} \text{ by first-order optimality}$$

$$\geq -\frac{\alpha(1-c_0)}{2}\|x_{t+1}(s) - \widehat{x}\|_2^2 + \frac{\alpha}{2}\|x_{t+1}(s) - \widehat{x}\|_2^2 - \frac{\alpha \varepsilon_\mathsf{g}^2}{2}$$

$$= \frac{c_0 \alpha}{2}\|x_{t+1}(s) - \widehat{x}\|_2^2 - \frac{\alpha \varepsilon_\mathsf{g}^2}{2}, \tag{40}$$

where the next-to-last step applies the initialization condition (19) (plus the fact that $\phi \geq 1$) to bound $\|\nabla \mathsf{g}(\widehat{x})\|^*$. On the other hand, we have

$$\mathsf{g}(x_{t+1}(s)) - \mathsf{g}(\widehat{x}) = \mathsf{g}(x_{t+1}(s)) - \mathsf{g}(x_t) + \mathsf{g}(x_t) - \mathsf{g}(\widehat{x})$$

$$\leq \langle x_{t+1}(s) - x_t, \nabla \mathsf{g}(x_t) \rangle + \frac{\beta}{2}\|x_{t+1}(s) - x_t\|_2^2 + \frac{\alpha \varepsilon_\mathsf{g}^2}{2} + \langle x_t - \widehat{x}, \nabla \mathsf{g}(x_t) \rangle - \frac{\alpha}{2}\|x_t - \widehat{x}\|_2^2 + \frac{\alpha \varepsilon_\mathsf{g}^2}{2}$$

$$= \langle x_{t+1}(s) - \widehat{x}, \nabla \mathsf{g}(x_t) \rangle + \frac{\beta}{2}\|x_{t+1}(s) - x_t\|_2^2 - \frac{\alpha}{2}\|x_t - \widehat{x}\|_2^2 + \alpha \varepsilon_\mathsf{g}^2, \tag{41}$$

where the inequality applies restricted strong convexity (11) and restricted smoothness (12). To bound the remaining inner product term, we have

$$\langle x_{t+1}(s) - \widehat{x}, \nabla \mathsf{g}(x_t) \rangle = \frac{1}{s}\langle x_{t+1}(s) - \widehat{x}, x_t - x'_{t+1}(s) \rangle$$

$$= \frac{1}{s}\langle x_{t+1}(s) - \widehat{x}, x_t - x_{t+1}(s) \rangle + \frac{1}{s}\langle x_{t+1}(s) - \widehat{x}, x_{t+1}(s) - x'_{t+1}(s) \rangle$$

$$\leq \frac{1}{s}\langle x_{t+1}(s) - \widehat{x}, x_t - x_{t+1}(s) \rangle + \frac{\alpha(1-c_0)}{2(\gamma^c + \gamma^d)}\left(\gamma^c \|\widehat{x} - x_{t+1}(s)\|_2^2 + \gamma^d \|\widehat{x} - x_t\|_2^2\right), \tag{42}$$

where the last step applies (39). For the first term on the right-hand side, we can trivially check that

$$\frac{1}{s}\langle x_{t+1}(s) - \widehat{x}, x_t - x_{t+1}(s) \rangle = \frac{1}{2s}\|x_t - \widehat{x}\|_2^2 - \frac{1}{2s}\|x_{t+1}(s) - \widehat{x}\|_2^2 - \frac{1}{2s}\|x_{t+1}(s) - x_t\|_2^2. \tag{43}$$

Combining steps (40), (41), (42), and (43), then, since $\frac{1}{2s} \geq \frac{1}{2\eta} = \frac{\beta}{2}$,

$$\frac{c_0 \alpha}{2}\|x_{t+1}(s) - \widehat{x}\|_2^2 \leq \frac{1}{2s}\|x_t - \widehat{x}\|_2^2 - \frac{1}{2s}\|x_{t+1}(s) - \widehat{x}\|_2^2$$

$$+ \frac{\alpha(1-c_0)}{2(\gamma^c + \gamma^d)}\left(\gamma^c \|\widehat{x} - x_{t+1}(s)\|_2^2 + \gamma^d \|\widehat{x} - x_t\|_2^2\right) - \frac{\alpha}{2}\|x_t - \widehat{x}\|_2^2 + 1.5\alpha \varepsilon_\mathsf{g}^2.$$

Rearranging terms we obtain

$$\|x_{t+1}(s) - \widehat{x}\|_2^2 \leq \left(1 - \frac{2\alpha c_0}{\frac{1}{s} + \alpha\left(c_0 - (1-c_0) \cdot \frac{\gamma^c}{\gamma^c + \gamma^d}\right)}\right)\|x_t - \widehat{x}\|_2^2 + \frac{3\alpha}{\frac{1}{s} + \alpha\left(c_0 - (1-c_0) \cdot \frac{\gamma^c}{\gamma^c + \gamma^d}\right)}\varepsilon_\mathsf{g}^2. \tag{44}$$

In particular, since $\|x_t - \widehat{x}\|_2 \leq \rho_0$ and $\varepsilon_g^2 \leq \frac{c_0 \rho_0^2}{1.5}$ by assumption, this proves that

$$\|x_{t+1}(s) - \widehat{x}\|_2 \leq \rho_0. \tag{45}$$

Therefore, we see that $s \in \mathcal{S}_0$.

To summarize, we have proved that for all $s \in \mathcal{S}$, we also have $s \in \mathcal{S}_0$; while for $s \in \mathcal{S}_0$, we also have $\min\{s + \Delta, \eta\} \in \mathcal{S}$, where $\Delta > 0$ is fixed. Starting with $s = 0 \in \mathcal{S}$ and proceeding inductively, we see that $\Delta \in \mathcal{S}$, then $2\Delta \in \mathcal{S}$, etc, until inductively we obtain $\eta \in \mathcal{S}$. Therefore, setting $s = \eta$ and $x_{t+1}(s) = x_{t+1}$, the above bounds (44) and (45) will hold. Looking at (44) in particular, since $s = \eta = 1/\beta$, we can simplify to

$$\|x_{t+1} - \widehat{x}\|_2^2 \leq \left(1 - \frac{2\alpha c_0}{\beta + \alpha\left(c_0 - (1-c_0) \cdot \frac{\gamma^c}{\gamma^c + \gamma^d}\right)}\right)\|x_t - \widehat{x}\|_2^2 + \frac{3\alpha}{\beta + \alpha\left(c_0 - (1-c_0) \cdot \frac{\gamma^c}{\gamma^c + \gamma^d}\right)}\varepsilon_g^2.$$

This proves that the inductive step (36) holds for $x_{t+1}$, as desired, which completes the proof of Theorem 4.

## B.3 Proof of necessity of the initialization conditions (Lemma 5)

By definition of the local concavity coefficients (6), since $x \notin \mathcal{C}_{\mathsf{dgn}}$ and $\gamma_x(\mathcal{C}) > 0$, we see that there must exist some $y \in \mathcal{C}$ and $z \in \mathbb{R}^d$, with $x = P_\mathcal{C}(z)$ such that

$$\langle y - x, z - x \rangle > \frac{\gamma_x(\mathcal{C})}{1+\epsilon}\|z - x\|^* \|x - y\|_2^2.$$

(If no such $y, z$ exist then local concavity coefficient at $x$ would be $\leq \frac{\gamma_x(\mathcal{C})}{1+\epsilon}$, which is a contradiction.) Next define $\mathsf{g}(v) = \frac{\alpha}{2}\|v - \widetilde{z}\|_2^2$, where $\widetilde{z} = x + \frac{z-x}{\|z-x\|^*} \cdot \frac{1+\epsilon}{2\gamma_x(\mathcal{C})}$. Then $\mathsf{g}$ is $\alpha$-strongly convex, and

$$2\gamma_x(\mathcal{C}) \cdot \|\nabla \mathsf{g}(x)\|^* = 2\gamma_x(\mathcal{C}) \cdot \|\alpha(\widetilde{z} - x)\|^* = \alpha(1+\epsilon).$$

Furthermore, for sufficiently small step size $\eta > 0$, we can see that $x - \eta \nabla \mathsf{g}(x)$ is a convex combination of $x$ and $z$, and so we have $P_\mathcal{C}(x - \eta\nabla\mathsf{g}(x)) = x$ by (29). Thus $x$ is a stationary point of projected gradient descent by (29), for sufficiently small $\eta > 0$. However, $x$ does not minimize $\mathsf{g}$ over $\mathcal{C}$, since

$$\|y - \widetilde{z}\|_2^2 - \|x - \widetilde{z}\|_2^2 = \|y - x\|_2^2 - 2\langle y - x, \widetilde{z} - x\rangle = \|y - x\|_2^2 - 2\langle y - x, z - x\rangle \cdot \frac{1+\epsilon}{2\gamma_x(\mathcal{C})\|z-x\|^*}$$

$$< \|y - x\|_2^2 - 2 \cdot \frac{\gamma_x(\mathcal{C})}{1+\epsilon}\|z - x\|^*\|x - y\|_2^2 \cdot \frac{1+\epsilon}{2\gamma_x(\mathcal{C})\|z-x\|^*} = 0,$$

which shows that $\mathsf{g}(y) < \mathsf{g}(x)$.

## B.4 Proof of equivalence of the exact and approximate settings (Lemma 6)

### B.4.1 Local concavity coefficients $\Rightarrow$ Family of approximate projections

Suppose that $\mathcal{C}$ has local concavity coefficients $\gamma_x(\mathcal{C})$ for all $x \in \mathcal{C}$. Then the inner product condition (5) for exact projection $P_\mathcal{C}$ gives

$$\langle y - u, z - u\rangle \leq \gamma_u(\mathcal{C})\|z - u\|^*\|y - u\|_2^2$$

for all $u, y \in \mathcal{C}$ and $z \in \mathbb{R}^d$ with $P_\mathcal{C}(z) = u$. Therefore, the general condition (21) holds for $P_x = P_\mathcal{C}$ when we set $\gamma_u^c(\mathcal{C}) = \gamma_u(\mathcal{C})$ and $\gamma_u^d(\mathcal{C}) = 0$.

### B.4.2 Family of approximate projections ⇒ Local concavity coefficients

Suppose there exists a family of operators $P_x : \mathbb{R}^d \to \mathcal{C}$ indexed over $x \in \mathcal{C}$, satisfying the inner product condition in (21). Assume that the local continuity property (22) holds. We now check that the local concavity coefficients satisfy $\gamma_x(\mathcal{C}) \leq \gamma_x^c(\mathcal{C}) + \gamma_x^d(\mathcal{C})$.

Fix any $x \in \mathcal{C}$. Assume that $\gamma_x^c(\mathcal{C}) + \gamma_x^d(\mathcal{C}) < \infty$ (otherwise, there is nothing to prove). We will verify that the inner product condition (5) holds at this point $x$ with $\gamma = \gamma_x^c(\mathcal{C}) + \gamma_x^d(\mathcal{C})$. Fix any $\epsilon > 0$, $y \in \mathcal{C}$, and $z \in \mathbb{R}^d$ such that $x = P_\mathcal{C}(z) \in \mathcal{C}$. Define
$$z_t = (1-t)x + tz.$$

We will take limits as $t$ approaches zero throughout the proof. Below we will prove that there exists some $t_0 > 0$ such that
$$P_x(z_t) = x \text{ for all } t \leq t_0. \tag{46}$$

Assume for now that this holds.

Take any $t \in (0, t_0)$. By the inner product condition (21) for the approximate projections $P_x$, since $x = P_x(z_t)$ for this small $t$,

$$\langle y - x, z_t - x \rangle = \langle y - P_x(z_t), z_t - P_x(z_t) \rangle$$
$$\leq \max\{\|z_t - P_x(z_t)\|^*, \|z_t - x\|^*\} \cdot \left( \gamma_{P_x(z_t)}^c(\mathcal{C}) \|y - P_x(z_t)\|_2^2 + \gamma_{P_x(z_t)}^d(\mathcal{C}) \|y - x\|_2^2 \right)$$
$$= \|z_t - x\|^* \cdot \left( \gamma_x^c(\mathcal{C}) + \gamma_x^d(\mathcal{C}) \right) \cdot \|y - x\|_2^2.$$

Plugging in the definition of $z_t$, and then dividing by $t$,
$$\langle y - x, z - x \rangle \leq \|z - x\|^* \cdot (\gamma_x^c(\mathcal{C}) + \gamma_x^d(\mathcal{C})) \|y - x\|_2^2,$$

proving that the inner product condition (5) holds, at this point $x$ and for any $y \in \mathcal{C}$, with $\gamma = \gamma_x^c(\mathcal{C}) + \gamma_x^d(\mathcal{C})$. By definition of the local concavity coefficients, if $x \in \mathcal{C}\backslash\mathcal{C}_{\text{dgn}}$, then $\gamma_x(\mathcal{C}) \leq \gamma = \gamma_x^c(\mathcal{C}) + \gamma_x^d(\mathcal{C})$. Now we consider the case that $x \in \mathcal{C}_{\text{dgn}}$. Since $x \mapsto \gamma_x^c(\mathcal{C}) + \gamma_x^d(\mathcal{C})$ is upper semi-continuous by assumption, we can find some $r > 0$ such that $\gamma := \sup_{x' \in \mathcal{C} \cap \mathbb{B}_2(x,r)} \left( \gamma_{x'}^c(\mathcal{C}) + \gamma_{x'}^d(\mathcal{C}) \right) < \infty$. By the reasoning above, then, the inner product condition (5) holds, with any $x' \in \mathcal{C} \cap \mathbb{B}_2(x,r)$ in place of $x$ and for any $y \in \mathcal{C}$, with the constant $\gamma$. However, if $x \in \mathcal{C}_{\text{dgn}}$ then this is not possible, due to Lemma 3. Thus we have reached a contradiction—if $\gamma_x^c(\mathcal{C}) + \gamma_x^d(\mathcal{C}) < \infty$ then we must have $x \in \mathcal{C}\backslash\mathcal{C}_{\text{dgn}}$. This completes the proof that $\gamma_x(\mathcal{C}) \leq \gamma_x^c(\mathcal{C}) + \gamma_x^d(\mathcal{C})$ for all $x \in \mathcal{C}$, assuming that (46) holds

Now it remains to be shown that (46) is indeed true. Since $z_t \to x$ trivially, and $P_x(z_t) \to x$ by the local continuity assumption (22), we see that $(z_t - P_x(z_t)) \to 0$. Since $u \mapsto \gamma_u^c(\mathcal{C})$ is upper semi-continuous, then $\gamma_{P_x(z_t)}^c \leq \gamma_x^c(\mathcal{C}) + \epsilon$ for sufficiently small $t > 0$, and also $\gamma_x^c(\mathcal{C}) < \infty$ by assumption. Therefore, for sufficiently small $t > 0$, we have $2\gamma_{P_x(z_t)}^c(\mathcal{C}) \max\{\|z_t - P_x(z_t)\|^*, \|z_t - x\|^*\} \leq 1$. Then, applying the inner product condition (21), with $u = P_x(z_t)$, with $y = x$, and with $z_t$ in place of $z$, we obtain

$$\|z_t - P_x(z_t)\|_2^2 - \|z_t - x\|_2^2 = 2\langle z_t - P_x(z_t), x - P_x(z_t) \rangle - \|x - P_x(z_t)\|_2^2$$
$$\leq 2\max\{\|z_t - P_x(z_t)\|^*, \|z_t - x\|^*\} \left( \gamma_{P_x(z_t)}^c(\mathcal{C}) \|x - P_x(z_t)\|_2^2 + \gamma_{P_x(z_t)}^d(\mathcal{C}) \|x - x\|_2^2 \right) - \|x - P_x(z_t)\|_2^2$$
$$= \|x - P_x(z_t)\|_2^2 \cdot (2\gamma_{P_x(z_t)}^c(\mathcal{C}) \max\{\|z_t - P_x(z_t)\|^*, \|z_t - x\|^*\} - 1) \leq 0,$$

implying that $\|z_t - P_x(z_t)\|_2^2 \leq \|z_t - x\|_2^2$. So,
$$0 \leq \|z_t - x\|_2^2 - \|z_t - P_x(z_t)\|_2^2 = 2\langle z_t - x, P_x(z_t) - x \rangle - \|P_x(z_t) - x\|_2^2 = 2t\langle z - x, P_x(z_t) - x \rangle - \|P_x(z_t) - x\|_2^2.$$

Furthermore, $x = P_\mathcal{C}(z)$, so
$$\|z - x\|_2^2 \leq \|z - P_x(z_t)\|_2^2 = \|z - x\|_2^2 + 2\langle z - x, x - P_x(z_t) \rangle + \|x - P_x(z_t)\|_2^2$$

so we have $2\langle z - x, P_x(z_t) - x \rangle \leq \|x - P_x(z_t)\|_2^2$. Since we have taken $t > 0$, we then have $\|x - P_x(z_t)\|_2^2 \leq t\|x - P_x(z_t)\|_2^2$ which implies that $x = P_x(z_t)$ for sufficiently small $t$. This proves (46), thus completing the proof of the lemma.

# C  Proofs for examples

In this section we prove results calculating the local concavity coefficients $\gamma_x(\mathcal{C})$ and the norm compatibility constant $\phi$ for the constraint sets considered in Section 5.

## C.1  Low rank constraints

For the low-rank constraint, we first recall the various matrix norms used in our analysis: the Frobenius norm $\|X\|_\mathsf{F} = \sqrt{\sum_{ij} X_{ij}^2}$ is the Euclidean $\ell_2$ norm when $X$ is reshaped into a vector; the nuclear norm $\|X\|_{\mathrm{nuc}} = \sum_i \sigma_i(X)$ is the sum of its singular values, and promotes a low-rank solution $X$ (in the same way that, for sparse vector estimation, the $\ell_1$ norm promotes sparse solutions); and the spectral norm $\|X\|_{\mathrm{sp}} = \sigma_1(X)$ is the largest singular value of $X$ (sometimes called the operator norm).

Recalling the subspace $T_X$ defined in (23) for any rank-$r$ matrix $X$, we begin with an auxiliary lemma:

**Lemma 14.** *Let $X, Y \in \mathbb{R}^{n \times m}$ satisfy $\mathrm{rank}(X), \mathrm{rank}(Y) \leq r$. Then*

$$\|P_{T_X}^\perp(Y)\|_{\mathrm{nuc}} \leq \frac{1}{2\sigma_r(X)} \|X - Y\|_\mathsf{F}^2.$$

*Proof of Lemma 14.* Assume $\sigma_r(X) > 0$ (otherwise the statement is trivial). For any matrix $M \in (T_X)^\perp$ with $\|M\|_{\mathrm{sp}} \leq 1$, define a function

$$\mathsf{f}_M(Z) = \frac{1}{2\sigma_r(X)} \|Z - X\|_\mathsf{F}^2 - \langle Z, M \rangle$$

over matrices $Z \in \mathbb{R}^{n \times m}$. We can rewrite this as

$$\mathsf{f}_M(Z) = \frac{1}{2\sigma_r(X)} \|Z - (X + \sigma_r(X) M)\|_\mathsf{F}^2 + \langle X, M \rangle - \frac{\sigma_r(X)}{2} \|M\|_\mathsf{F}^2.$$

Now, we minimize $\mathsf{f}_M(Z)$ over a rank constraint:

$$\underset{\mathrm{rank}(Z) \leq r}{\arg\min}\, \mathsf{f}_M(Z) = \arg\min_Z \{\|Z - (X + \sigma_r(X) M)\|_\mathsf{F}^2 : \mathrm{rank}(Z) \leq r\} = P_\mathcal{C}\left(X + \sigma_r(X) M\right).$$

Since $\sigma_1(X), \ldots, \sigma_r(X) \geq \sigma_r(X)$ while $\|\sigma_r(X) M\|_{\mathrm{sp}} \leq \sigma_r(X)$, and $M \in (T_X)^\perp$, we see that

$$X = P_\mathcal{C}\left(X + \sigma_r(X) M\right).$$

(It may be the case that $X$ and $\sigma_r(X) M$ both have one or more singular values exactly equal to $\sigma_r(X)$, in which case the projection is not unique, but $X$ is always one of the values of the projection.) So, $Z = X$ minimizes $\mathsf{f}_M(Z)$ over rank-$r$ matrices, and therefore, for any $Z$ with $\mathrm{rank}(Z) \leq r$,

$$\mathsf{f}_M(Z) \geq \mathsf{f}_M(X) = \frac{1}{2\sigma_r(X)} \|X - X\|_\mathsf{F}^2 + \langle X, M \rangle = 0,$$

since $\langle X, M \rangle = 0$ due to $M \in (T_X)^\perp$. Therefore, in particular, $\mathsf{f}_M(Y) \geq 0$ which implies that

$$\langle Y, M \rangle \leq \frac{1}{2\sigma_r(X)} \|Y - X\|_\mathsf{F}^2.$$

Since this is true for all $M \in (T_X)^\perp$ with $\|M\|_{\mathrm{sp}} \leq 1$, we have proved that

$$\|P_{T_X}^\perp(Y)\|_{\mathrm{nuc}} = \max_{M \in (T_X)^\perp, \|M\|_{\mathrm{sp}} = 1} \langle Y, M \rangle \leq \frac{1}{2\sigma_r(X)} \|Y - X\|_\mathsf{F}^2,$$

as desired. □

*Proof of Lemma 7.* First, let $P_\mathcal{C}(Z)$ be any closest rank-$r$ matrix to $Z$ (not necessarily unique), and let $U \in \mathbb{R}^{n \times r}$ and $V \in \mathbb{R}^{m \times r}$ be orthonormal bases for the column span and row span of $P_\mathcal{C}(Z)$ (that is, if $P_\mathcal{C}(Z)$ is unique then the columns of $U$ and $V$ are the top $r$ left and right singular vectors of $Z$). Regardless of uniqueness we will have $Z - P_\mathcal{C}(Z)$ orthogonal to $U$ on the left and to $V$ on the right, i.e. we can write

$$Z - P_\mathcal{C}(Z) = (\mathbf{I} - UU^\top) \cdot (Z - P_\mathcal{C}(Z)) \cdot (\mathbf{I} - VV^\top).$$

We then have

$$\begin{aligned}
\langle Y - P_\mathcal{C}(Z), Z - P_\mathcal{C}(Z) \rangle &= \langle Y - P_\mathcal{C}(Z), (\mathbf{I} - UU^\top) \cdot (Z - P_\mathcal{C}(Z)) \cdot (\mathbf{I} - VV^\top) \rangle \\
&= \langle (\mathbf{I} - UU^\top) \cdot (Y - P_\mathcal{C}(Z)) \cdot (\mathbf{I} - VV^\top), Z - P_\mathcal{C}(Z) \rangle \\
&\leq \|(\mathbf{I} - UU^\top) \cdot (Y - P_\mathcal{C}(Z)) \cdot (\mathbf{I} - VV^\top)\|_{\text{nuc}} \cdot \|Z - P_\mathcal{C}(Z)\|_{\text{sp}} \\
&\leq \|(\mathbf{I} - UU^\top) \cdot Y \cdot (\mathbf{I} - VV^\top)\|_{\text{nuc}} \cdot \|Z - P_\mathcal{C}(Z)\|_{\text{sp}},
\end{aligned}$$

where the last step holds since $P_\mathcal{C}(Z)$ is spanned by $U$ on the left and $V$ on the right. Applying Lemma 14 with $X = P_\mathcal{C}(Z)$, which trivially has $U, V$ as its left and right singular vectors, we obtain

$$\|(\mathbf{I} - UU^\top) \cdot Y \cdot (\mathbf{I} - VV^\top)\|_{\text{nuc}} \leq \frac{1}{2\sigma_r(P_\mathcal{C}(Z))} \|Y - P_\mathcal{C}(Z)\|_{\mathsf{F}}^2.$$

Therefore,

$$\langle Y - P_\mathcal{C}(Z), Z - P_\mathcal{C}(Z) \rangle \leq \frac{1}{2\sigma_r(P_\mathcal{C}(Z))} \|Z - P_\mathcal{C}(Z)\|_{\text{sp}} \|Y - P_\mathcal{C}(Z)\|_{\mathsf{F}}^2.$$

This proves that $\gamma_X(\mathcal{C}) \leq \frac{1}{2\sigma_r(X)}$ for all $x \in \mathcal{C}$ by the inner product condition (5).

To prove equality, take any $X \in \mathcal{C} \backslash \mathcal{C}_{\mathsf{dgn}}$ (that is, we assume that $\text{rank}(X) = r$), and let $X = \sigma_1 u_1 v_1^\top + \cdots + \sigma_r u_r v_r^\top$ be a singular value decomposition with $\sigma_1 \geq \cdots \geq \sigma_r > 0$. Let $u' \in \mathbb{R}^n, v' \in \mathbb{R}^m$ be unit vectors orthogonal to the left and right singular vectors of $X$, respectively. Define

$$Y = \sigma_1 u_1 v_1^\top + \cdots + \sigma_{r-1} u_{r-1} v_{r-1}^\top + \sigma_r u' {v'}^\top$$

and

$$Z = X + c u' {v'}^\top,$$

for some fixed $c \in (0, \sigma_r)$. Then $P_\mathcal{C}(Z) = X$, and we have

$$\langle Y - X, Z - X \rangle = \langle \sigma_r u' {v'}^\top - \sigma_r u_r v_r^\top, c u' {v'}^\top \rangle = c \sigma_r$$

while

$$\|Z - X\|_{\text{sp}} \|Y - X\|_{\mathsf{F}}^2 = \|c u' {v'}^\top\|_{\text{sp}} \|\sigma_r u' {v'}^\top - \sigma_r u_r v_r^\top\|_{\mathsf{F}}^2 = 2 c \sigma_r^2,$$

therefore by the inner product condition (5), we must have $\gamma_X(\mathcal{C}) \geq \frac{1}{2\sigma_r(X)}$.

Turning to the norm compatibility condition, the desired bound is an immediate result of the Eckart–Young theorem [11], as

$$\|Z - P_\mathcal{C}(Z)\|_{\text{sp}} \leq \|Z - X\|_{\text{sp}},$$

for all $X \in \mathcal{C}$ and $Z \in \mathbb{R}^{n \times m}$. □

*Proof of Lemma 8.* Let $X$ be any matrix with $\text{rank}(X) \leq r$ and let $Z$ be any matrix. Assume $X, P_X(Z) \in \mathbb{B}_2(\widehat{X}, \rho)$ for $\rho = \frac{\sigma_r(\widehat{X})}{4}$. According to Weyl's inequality, we will have $\sigma_r(X), \sigma_r(P_X(Z)) \geq \frac{3\sigma_r(\widehat{X})}{4}$. Write $T = T_X$ for convenience, and define $Z_T = P_T(Z)$ and $Z_\perp = P_T^\perp(Z)$.

Then $P_X(Z) = P_\mathcal{C}(P_T(Z)) = P_\mathcal{C}(Z_T)$, and so

$$\begin{aligned}
\langle \widehat{X} - P_X(Z), Z - P_X(Z) \rangle &= \langle \widehat{X} - P_\mathcal{C}(Z_T), Z - P_\mathcal{C}(Z_T) \rangle \\
&= \underbrace{\langle \widehat{X} - P_\mathcal{C}(Z_T), P_T(Z - P_\mathcal{C}(Z_T)) \rangle}_{\text{(Term 1)}} + \underbrace{\langle \widehat{X} - P_\mathcal{C}(Z_T), P_T^\perp(Z - P_\mathcal{C}(Z_T)) \rangle}_{\text{(Term 2)}}.
\end{aligned}$$

First consider (Term 1). We have

$$\langle \widehat{X} - P_\mathcal{C}(Z_T), P_T(Z - P_\mathcal{C}(Z_T))\rangle$$
$$= \langle \widehat{X} - P_\mathcal{C}(Z_T), Z_T - P_T(P_\mathcal{C}(Z_T))\rangle$$
$$= \langle \widehat{X} - P_\mathcal{C}(Z_T), Z_T - P_\mathcal{C}(Z_T)\rangle + \langle \widehat{X} - P_\mathcal{C}(Z_T), P_T^\perp(P_\mathcal{C}(Z_T))\rangle$$
$$= \langle \widehat{X} - P_\mathcal{C}(Z_T), Z_T - P_\mathcal{C}(Z_T)\rangle - \langle \widehat{X} - P_\mathcal{C}(Z_T), P_T^\perp(Z_T - P_\mathcal{C}(Z_T))\rangle$$
$$\leq \frac{2}{3\sigma_r(\widehat{X})}\|Z_T - P_\mathcal{C}(Z_T)\|_{\mathrm{sp}}\|\widehat{X} - P_\mathcal{C}(Z_T)\|_{\mathsf{F}}^2 - \langle \widehat{X} - P_\mathcal{C}(Z_T), P_T^\perp(Z_T - P_\mathcal{C}(Z_T))\rangle$$
$$\leq \frac{2}{3\sigma_r(\widehat{X})}\|Z_T - P_\mathcal{C}(Z_T)\|_{\mathrm{sp}}\|\widehat{X} - P_\mathcal{C}(Z_T)\|_{\mathsf{F}}^2 + \|Z_T - P_\mathcal{C}(Z_T)\|_{\mathrm{sp}}\left(\|P_T^\perp(\widehat{X})\|_{\mathrm{nuc}} + \|P_T^\perp(P_\mathcal{C}(Z_T))\|_{\mathrm{nuc}}\right)$$
$$\leq \frac{2}{3\sigma_r(\widehat{X})}\|Z_T - P_\mathcal{C}(Z_T)\|_{\mathrm{sp}}\|\widehat{X} - P_\mathcal{C}(Z_T)\|_{\mathsf{F}}^2 + \|Z_T - P_\mathcal{C}(Z_T)\|_{\mathrm{sp}}\left(\frac{2}{3\sigma_r(\widehat{X})}\|\widehat{X} - X\|_{\mathsf{F}}^2 + \frac{2}{3\sigma_r(\widehat{X})}\|P_\mathcal{C}(Z_T) - X\|_{\mathsf{F}}^2\right),$$

where the first inequality applies the inner product condition (5), using the fact that $\gamma_{P_\mathcal{C}(Z_T)} = \frac{1}{2\sigma_r(P_\mathcal{C}(Z_T))} \leq \frac{2}{3\sigma_r(\widehat{X})}$; the second inequality uses the duality between nuclear norm and spectral norm; and the third applies Lemma 14 to both nuclear norm terms since $\mathrm{rank}(\widehat{X}), \mathrm{rank}(P_\mathcal{C}(Z_T)) \leq r$ and $\frac{1}{2\sigma_r(X)} \leq \frac{2}{3\sigma_r(\widehat{X})}$. Also, since

$$\|P_\mathcal{C}(Z_T) - X\|_{\mathsf{F}}^2 \leq 2\|P_\mathcal{C}(Z_T) - \widehat{X}\|_{\mathsf{F}}^2 + 2\|\widehat{X} - X\|_{\mathsf{F}}^2,$$

we can simplify our bound to

$$\langle \widehat{X} - P_\mathcal{C}(Z_T), P_T(Z - P_\mathcal{C}(Z_T))\rangle \leq \frac{2}{\sigma_r(\widehat{X})}\|Z_T - P_\mathcal{C}(Z_T)\|_{\mathrm{sp}}\left(\|\widehat{X} - P_\mathcal{C}(Z_T)\|_{\mathsf{F}}^2 + \|\widehat{X} - X\|_{\mathsf{F}}^2\right).$$

Finally, we have

$$\|Z_T - P_\mathcal{C}(Z_T)\|_{\mathrm{sp}} \leq \|Z - P_\mathcal{C}(Z_T)\|_{\mathrm{sp}} + \|P_T^\perp(Z)\|_{\mathrm{sp}} = \|Z - P_X(Z)\|_{\mathrm{sp}} + \|P_T^\perp(Z - X)\|_{\mathrm{sp}}$$
$$\leq \|Z - P_X(Z)\|_{\mathrm{sp}} + \|Z - X\|_{\mathrm{sp}} \leq 2\max\{\|Z - P_X(Z)\|_{\mathrm{sp}}, \|Z - X\|_{\mathrm{sp}}\}$$

since $X \in T$ by definition and $P_T^\perp$ is contractive with respect to spectral norm. Then, returning to the work above,

$$\langle \widehat{X} - P_\mathcal{C}(Z_T), P_T(Z - P_\mathcal{C}(Z_T))\rangle \leq \frac{4}{\sigma_r(\widehat{X})}\max\{\|Z - P_X(Z)\|_{\mathrm{sp}}, \|Z - X\|_{\mathrm{sp}}\}\left(\|\widehat{X} - P_\mathcal{C}(Z_T)\|_{\mathsf{F}}^2 + \|\widehat{X} - X\|_{\mathsf{F}}^2\right).$$

Next we turn to (Term 2). We have

$$\langle \widehat{X} - P_\mathcal{C}(Z_T), P_T^\perp(Z - P_\mathcal{C}(Z_T))\rangle = \langle P_T^\perp(\widehat{X}), Z - P_\mathcal{C}(Z_T)\rangle - \langle P_T^\perp(P_\mathcal{C}(Z_T)), Z - P_\mathcal{C}(Z_T)\rangle$$
$$\leq \|P_T^\perp(\widehat{X})\|_{\mathrm{nuc}}\|Z - P_\mathcal{C}(Z_T)\|_{\mathrm{sp}} + \|P_T^\perp(P_\mathcal{C}(Z_T))\|_{\mathrm{nuc}}\|Z - P_\mathcal{C}(Z_T)\|_{\mathrm{sp}}$$
$$\leq \frac{2}{3\sigma_r(\widehat{X})}\|\widehat{X} - X\|_{\mathsf{F}}^2\|Z - P_\mathcal{C}(Z_T)\|_{\mathrm{sp}} + \frac{2}{3\sigma_r(\widehat{X})}\|P_\mathcal{C}(Z_T) - X\|_{\mathsf{F}}^2\|Z - P_\mathcal{C}(Z_T)\|_{\mathrm{sp}}$$
$$\leq \frac{2}{\sigma_r(\widehat{X})}\|Z - P_\mathcal{C}(Z_T)\|_{\mathrm{sp}}\left(\|\widehat{X} - X\|_{\mathsf{F}}^2 + \|\widehat{X} - P_\mathcal{C}(Z_T)\|_{\mathsf{F}}^2\right),$$

where the second inequality applies Lemma 14 as before, while the third inequality again uses

$$\|P_\mathcal{C}(Z_T) - X\|_{\mathsf{F}}^2 \leq 2\|P_\mathcal{C}(Z_T) - \widehat{X}\|_{\mathsf{F}}^2 + 2\|\widehat{X} - X\|_{\mathsf{F}}^2.$$

Putting the bounds for (Term 1) and (Term 2) together, we conclude that

$$\langle \widehat{X} - P_X(Z), Z - P_X(Z)\rangle \leq \frac{6}{\sigma_r(\widehat{X})}\cdot\max\{\|Z - P_X(Z)\|_{\mathrm{sp}}, \|Z - X\|_{\mathrm{sp}}\}\left(\|\widehat{X} - P_\mathcal{C}(Z_T)\|_{\mathsf{F}}^2 + \|\widehat{X} - X\|_{\mathsf{F}}^2\right),$$

thus proving the inner product condition (17).

Now we turn to the norm compatibility condition (18). We have

$$\|Z - P_X(Z)\|_{\text{sp}} = \|Z - P_{\mathcal{C}}(Z_T)\|_{\text{sp}} \leq \|Z_T - P_{\mathcal{C}}(Z_T)\|_{\text{sp}} + \|Z_\perp\|_{\text{sp}} \leq \|Z_T - X\|_{\text{sp}} + \|Z_\perp\|_{\text{sp}},$$

where the last step holds by the Eckart–Young theorem [11]. Next, since $Z_T = Z - Z_\perp$, we have

$$\|Z_T - X\|_{\text{sp}} \leq \|Z - X\|_{\text{sp}} + \|Z_\perp\|_{\text{sp}},$$

while since $P_T^\perp(X) = 0$, we have

$$\|Z_\perp\|_{\text{sp}} = \|P_T^\perp(Z - X)\|_{\text{sp}} \leq \|Z - X\|_{\text{sp}}.$$

Combining everything, then,

$$\|Z - P_X(Z)\|_{\text{sp}} \leq 3\|Z - X\|_{\text{sp}},$$

for any $X \in \mathcal{C}$ and any $Z \in \mathbb{R}^{n \times m}$. This proves that the norm compatibility condition holds with $\phi = 3$.

Finally, we consider the local continuity condition (20). Fix any $c, \epsilon > 0$ and any $X \in \mathcal{C}$ and $Z \in \mathbb{R}^{n \times m}$ so that $\|X - \widehat{X}\|_{\mathsf{F}} \leq \rho$ and $\|P_X(Z) - \widehat{X}\|_{\mathsf{F}} \leq \rho$ where again $\rho = \frac{\sigma_r(\widehat{X})}{4}$. Suppose that

$$2(\gamma^{\text{c}} + \gamma^{\text{d}})\|Z - P_X(Z)\|_{\text{sp}} = \frac{24}{\sigma_r(\widehat{X})}\|Z - P_X(Z)\|_{\text{sp}} \leq 1 - c$$

and take any $W \in \mathbb{R}^{n \times m}$ with $\|Z - W\|_{\mathsf{F}} \leq \delta := \epsilon/4.5$.

Then we calculate

$$\gamma_{P_X(Z)}(\mathcal{C}) \leq \frac{1}{2\sigma_r(P_X(Z))} \leq \frac{2}{3\sigma_r(\widehat{X})},$$

as before. And,

$$\|P_T(Z) - P_X(Z)\|_{\text{sp}} \leq \|Z - P_X(Z)\|_{\text{sp}} + \|P_T^\perp(Z)\|_{\text{sp}} \quad \text{by the triangle inequality}$$
$$= \|Z - P_X(Z)\|_{\text{sp}} + \|P_T^\perp(Z - X)\|_{\text{sp}} \quad \text{since } P_T^\perp(X) = 0$$
$$\leq \|Z - P_X(Z)\|_{\text{sp}} + \|Z - X\|_{\text{sp}} \quad \text{since } P_T^\perp \text{ is contractive with respect to spectral norm}$$
$$\leq 2\|Z - P_X(Z)\|_{\text{sp}} + \|P_X(Z) - \widehat{X}\|_{\text{sp}} + \|X - \widehat{X}\|_{\text{sp}} \quad \text{by the triangle inequality}$$
$$\leq 2 \cdot \frac{(1-c)\sigma_r(\widehat{X})}{24} + 2\rho \quad \text{since } \|\cdot\|_{\text{sp}} \leq \|\cdot\|_{\mathsf{F}}$$
$$\leq \frac{7\sigma_r(\widehat{X})}{12}.$$

Then

$$\|P_X(Z) - P_X(W)\|_{\mathsf{F}} = \|P_{\mathcal{C}}(P_T(Z)) - P_{\mathcal{C}}(P_T(W))\|_{\mathsf{F}}$$
$$\leq \frac{\|P_T(Z) - P_T(W)\|_{\mathsf{F}}}{1 - 2\gamma_{P_X(Z)}(\mathcal{C})\|P_T(Z) - P_{\mathcal{C}}(P_T(Z))\|_{\text{sp}}} \quad \text{by the contraction property (3)}$$
$$\leq \frac{\|P_T(Z) - P_T(W)\|_{\mathsf{F}}}{1 - 2 \cdot \frac{2}{3\sigma_r(\widehat{X})} \cdot \frac{7\sigma_r(\widehat{X})}{12}} \quad \text{by the calculations above}$$
$$= 4.5\|P_T(Z) - P_T(W)\|_{\mathsf{F}}$$
$$\leq 4.5\|Z - W\|_{\mathsf{F}},$$

since $T$ is a subspace so $P_T$ is contractive with respect to the Frobenius norm. Since $\|Z - W\|_{\mathsf{F}} \leq \delta = \epsilon/4.5$ by assumption, this proves that $\|P_X(Z) - P_X(W)\|_{\mathsf{F}} \leq \epsilon$, as desired. $\square$

## C.2 Sparsity

*Proof of Lemma 9.* We check the local concavity coefficients. Fix any $x \in \mathcal{C}$. As before, if $x$ is in the interior (i.e. $\text{Pen}(x) < c$) then $\gamma_x(\mathcal{C}) = 0$, so we turn to the case that $\text{Pen}(x) = c$, and in particular, $x \neq 0$. Without loss of generality, assume that $x_1 > 0$ and that $x_1$ is the smallest nonzero coordinate of $x$ (and then $x_{\min} = x_1$). Choose any $y \in \mathcal{C}$ and $t \in [0, 1]$. Let

$$x_t = (1-t)x + ty \text{ and } z_t = x_t - s_t \mathbf{e}_1$$

where $\mathbf{e}_1 = (1, 0, \ldots, 0)$ and

$$s_t = t \cdot \frac{\mu/2}{\mathsf{p}'((x_t)_1)} \cdot \|x - y\|_2^2.$$

Since $\lim_{t \searrow} x_t = x$, and $\mathsf{p}$ is continuously differentiable (since it is both concave and differentiable on the positive real line), we have

$$\lim_{t \searrow 0} \frac{s_t}{t} = \frac{\mu/2}{\mathsf{p}'(x_1)} \cdot \|x - y\|_2^2.$$

In particular, this implies that, for sufficiently small $t$, we have $(x_t)_1 > 0$ and $(z_t)_1 > 0$.

We claim that $\text{Pen}(z_t) \leq c$, in which case

$$\lim_{t \searrow 0} \frac{\min_{x' \in \mathcal{C}} \|x_t - x'\|_1}{t} \leq \lim_{t \searrow 0} \frac{\|x_t - z_t\|_1}{t} = \lim_{t \searrow 0} \frac{s_t}{t} = \frac{\mu/2}{\mathsf{p}'(x_{\min})} \cdot \|x - y\|_2^2,$$

which proves the lemma.

It now remains to check that $\text{Pen}(z_t) \leq c$. We have, for coordinate $i = 1$,

$$\mathsf{p}(|z_t|_i) = \mathsf{p}((x_t)_1 - s_t) \leq \mathsf{p}((x_t)_1) - s_t \mathsf{p}'((x_t)_1),$$

since $0 < (x_t)_1 - s_t < (x_t)_1$ and $\rho$ is concave over $\mathbb{R}_+$. And, for every coordinate $i$,

$$\begin{aligned}
\mathsf{p}(|(x_t)_i|) &= \mathsf{p}\left(|(1-t)x_i + ty_i|\right) \\
&\leq \mathsf{p}\left((1-t)|x_i| + t|y_i|\right) \quad \text{since } \rho \text{ is nondecreasing} \\
&\leq (1-t)\mathsf{p}(|x_i|) + t\mathsf{p}(|y_i|) + \frac{\mu}{2}t(1-t)(|x_i| - |y_i|)^2 \quad \text{since } t \mapsto \mathsf{p}(t) + \mu t^2/2 \text{ is convex} \\
&\leq (1-t)\mathsf{p}(|x_i|) + t\mathsf{p}(|y_i|) + \frac{\mu}{2}t(1-t)(x_i - y_i)^2.
\end{aligned}$$

Therefore,

$$\begin{aligned}
\text{Pen}(z_t) = \sum_i \mathsf{p}(|z_t|_i) &\leq \left(\sum_i (1-t)\mathsf{p}(|x_i|) + t\mathsf{p}(|y_i|) + \frac{\mu}{2}t(1-t)(x_i - y_i)^2\right) - s_t \mathsf{p}'((x_t)_1) \\
&\leq (1-t)\text{Pen}(x) + t\text{Pen}(y) + \frac{\mu}{2}\|x - y\|_2^2 - s_t \mathsf{p}'((x_t)_1) \leq c + \frac{\mu}{2}\|x - y\|_2^2 - s_t \mathsf{p}'((x_t)_1) = c,
\end{aligned}$$

where the last step holds by definition of $s_t$. $\square$

*Proof of Lemma 10.* Choose any $z \in \mathbb{R}^d$ and let $x = P_\mathcal{C}(z)$. Without loss of generality we take $z, x \geq 0$ and consider only the nontrivial case that $z \notin \mathcal{C}$, therefore $\text{Pen}(x) = c$ as $x$ cannot lie in the interior of the set. Furthermore, we can see that

$$x_i = (P_\mathcal{C}(z))_i = \max\{0, z_i - \lambda \mathsf{p}'(x_i)\} \text{ for all } i = 1, \ldots, d$$

for some $\lambda \geq 0$ by optimality of $x$ as the closest point to $z$ under the constraint $\sum_i \mathsf{p}(x_i) \leq c$. That is, $P_\mathcal{C}(z)$ behaves like projection to a weighted $\ell_1$ ball, with weights determined by the projection $x_i$ itself. Note that $\mathsf{Pen}(x) = c$, therefore $\mathsf{p}(x_i) \leq c$ and so $\mathsf{p}'(x_i) \geq \mathsf{p}'(\mathsf{p}^{-1}(c))$ for all $i$. Now consider any $w$ with $\|w - z\|_\infty < \lambda \mathsf{p}'(\mathsf{p}^{-1}(c))$. Let $v$ be the vector with entries $v_i = \min\{\max\{0, w_i\}, z_i\}$. Then for all $i \in S$, we have

$$|v_i - z_i| < \lambda \mathsf{p}'(\mathsf{p}^{-1}(c)) \leq \lambda \mathsf{p}'(x_i) = |x_i - z_i|.$$

And, for $i \notin S$, $|v_i - z_i| \leq |z_i| = |x_i - z_i|$. So, $\|v - z\|_2 < \|x - z\|_2$ (since $x \neq 0$ and so $S \neq \varnothing$). Therefore, since $x$ is the closest point to $z$ in $\mathcal{C}$, we must have $v \notin \mathcal{C}$. Since $|v_i| \leq |w_i|$ for all $i$ by construction of $v$, this means that $\mathsf{Pen}(w) \geq \mathsf{Pen}(v) > c$, and therefore any $w$ with $\|w - z\|_\infty < \lambda \mathsf{p}'(\mathsf{p}^{-1}(c))$ cannot lie in $\mathcal{C}$. In other words,

$$\min_{w \in \mathcal{C}} \|w - z\|_\infty \geq \lambda \mathsf{p}'(\mathsf{p}^{-1}(c)).$$

On the other hand,

$$\|z - x\|_\infty = \max_i |z_i - \max\{0, z_i - \lambda \mathsf{p}'(x_i)\}| \leq \lambda \max_i \mathsf{p}'(x_i) \leq \lambda,$$

since $\mathsf{p}'(t) \leq 1$ for all $t$. So, the norm compatibility condition (13) is satisfied with $\phi = \frac{1}{\mathsf{p}'(\mathsf{p}^{-1}(c))}$. $\square$

## C.3 Other examples

*Proof of Lemma 11.* Let $X, Y \in \mathcal{C}$. For a fixed $t \in (0, 1)$, let $(1 - t)X + tY = ADB^\top$ be a singular value decomposition. Since $AB^\top \in \mathbb{R}^{n \times r}$ is an orthonormal matrix, we then have

$$\min_{Z \in \mathcal{C}} \|Z - ((1-t)X + tY)\|_{\mathsf{nuc}} \leq \|AB^\top - ((1-t)X + tY)\|_{\mathsf{nuc}} = \|AB^\top - ADB^\top\|_{\mathsf{nuc}} = \sum_{i=1}^r (1 - D_{ii}).$$

Furthermore,

$$\begin{aligned}
\|D\|_{\mathsf{F}}^2 &= \|(1-t)X + tY\|_{\mathsf{F}}^2 \\
&= (1-t)^2 \|X\|_{\mathsf{F}}^2 + t^2 \|Y\|_{\mathsf{F}}^2 + 2t(1-t)\langle X, Y \rangle \\
&= (1-t)^2 \|X\|_{\mathsf{F}}^2 + t^2 \|Y\|_{\mathsf{F}}^2 + t(1-t)\left(\|X\|_{\mathsf{F}}^2 + \|Y\|_{\mathsf{F}}^2 - \|X - Y\|_{\mathsf{F}}^2\right) \\
&= (1-t)^2 r + t^2 r + t(1-t)\left(r + r - \|X - Y\|_{\mathsf{F}}^2\right) \\
&= r - t(1-t)\|X - Y\|_{\mathsf{F}}^2.
\end{aligned}$$

A trivial calculation shows that $1 - D_{ii} = \frac{1 - D_{ii}^2}{2} + \frac{(1 - D_{ii})^2}{2}$, so we have

$$\min_{Z \in \mathcal{C}} \|Z - ((1-t)X + tY)\|_{\mathsf{nuc}} \leq \sum_{i=1}^r (1 - D_{ii}) = \sum_{i=1}^r \frac{1 - D_{ii}^2}{2} + \frac{(1 - D_{ii})^2}{2} = \frac{r - \|D\|_{\mathsf{F}}^2}{2} + \sum_{i=1}^r \frac{(1 - D_{ii})^2}{2}$$

$$= \frac{1}{2} t(1-t) \|X - Y\|_{\mathsf{F}}^2 + \sum_{i=1}^r \frac{(1 - D_{ii})^2}{2}.$$

Furthermore, we can show that the last term is $o(t)$, as follows. For any unit vector $u \in \mathbb{R}^r$,

$$\|((1-t)X + tY)u\|_2 \geq (1-t)\|Xu\|_2 - t\|Yu\|_2 \geq 1 - 2t$$

since $X, Y$ are both orthonormal. Therefore $(1-t)X + tY$ has all its singular values $\geq 1 - 2t$, that is, $D_{ii} \geq 1 - 2t$ for all $i$. And trivially $\|((1-t)X + tY)u\|_2 \leq 1$ so $D_{ii} \leq 1$. Then $\sum_{i=1}^r (1 - D_{ii})^2 \leq \sum_{i=1}^r (2t)^2 = 4t^2 r$, so we have

$$\min_{Z \in \mathcal{C}} \|Z - ((1-t)X + tY)\|_{\mathsf{nuc}} \leq \frac{1}{2} t(1-t) \|X - Y\|_{\mathsf{F}}^2 + 2t^2 r.$$

Dividing by $t$ and taking a limit,

$$\lim_{t \searrow 0} \frac{\min_{Z \in \mathcal{C}} \|Z - ((1-t)X + tY)\|_{\mathsf{nuc}}}{t} \leq \frac{1}{2} \|X - Y\|_{\mathsf{F}}^2.$$

Comparing to the curvature condition (1) we see that $\gamma_X(\mathcal{C}) \leq \frac{1}{2}$, as desired.

Next, to obtain equality, take any $X \in \mathcal{C}$. Fix any $c \in (0,1)$. Let $Y = -X \in \mathcal{C}$ and $Z = cX \in \mathbb{R}^{n \times r}$. Clearly, $P_\mathcal{C}(Z) = X$. By the contraction property (3), we must have

$$(1 - \gamma_X(\mathcal{C})\|Z - X\|^*)\|Y - X\|_\mathsf{F} \leq \|Y - Z\|_\mathsf{F}.$$

Plugging in our choices for $Y$ and $Z$, we obtain

$$(1 - \gamma_X(\mathcal{C}) \cdot (1-c)) \cdot 2\sqrt{r} \leq (1+c)\sqrt{r},$$

and so $\gamma_X(\mathcal{C}) \geq \frac{1}{2}$.

Now we check the norm compatibility condition. For any $X \in \mathbb{R}^{n \times r}$, write $X = ADB^\top$. Then $P_\mathcal{C}(X) = AB^\top$ and so $\|X - P_\mathcal{C}(X)\|^* = \|ADB^\top - AB^\top\|_{\mathrm{sp}} = \max\{d_1 - 1, 1 - d_r\}$, where $d_1 \geq \cdots \geq d_r \geq 0$ are the diagonal entries of $D$, i.e. the singular values of $X$. Let $u \in \mathbb{R}^d$ be the first column of $B$, so that $\|Xu\|_2 = d_1$. Then for any $W \in \mathcal{C}$, we have $\|Wu\|_2 = 1$ since $W$ is orthonormal and $u$ is a unit vector, so

$$\|X - W\|^* = \|X - W\|_{\mathrm{sp}} \geq \|(X-W)u\|_2 \geq \|Xu\|_2 - \|Wu\|_2 = d_1 - 1.$$

Now let $v \in \mathbb{R}^d$ be the $r$th column of $B$, so that $\|Xv\|_2 = d_r$. Similarly we have

$$\|X - W\|^* = \|X - W\|_{\mathrm{sp}} \geq \|(X-W)v\|_2 \geq \|Wv\|_2 - \|Xv\|_2 = 1 - d_r.$$

Therefore, $\|X - W\|^* \geq \max\{d_1 - 1, 1 - d_r\} = \|X - P_\mathcal{C}(X)\|^*$, proving that $\phi = 1$. $\square$

*Proof of Lemma 12.* For $X, Y \in \mathcal{C}$, write $X = UU^\top$, and $Y = VV^\top$ for some orthonormal matrices $U, V \in \mathbb{R}^{n \times r}$. For $t \in (0,1)$, let $U_t = (1-t)U + tV$, and let $U_t = ADB^\top$ be a singular value decomposition. Then $AB^\top$ is the projection of $U_t$ onto the set of orthonormal $n \times r$ matrices. Since $A \in \mathbb{R}^{n \times r}$ is orthonormal, we have $AA^\top \in \mathcal{C}$, and so

$$\min_{Z \in \mathcal{C}} \|Z - ((1-t)X + tY)\|_{\mathrm{nuc}} \leq \|AA^\top - ((1-t)X + tY)\|_{\mathrm{nuc}}$$

$$\leq \underbrace{\|AA^\top - U_tU_t^\top\|_{\mathrm{nuc}}}_{\text{(Term 1)}} + \underbrace{\|U_tU_t^\top - ((1-t)X + tY)\|_{\mathrm{nuc}}}_{\text{(Term 2)}}.$$

For (Term 1),

$$\|AA^\top - U_tU_t^\top\|_{\mathrm{nuc}} = \|AA^\top - ADB^\top \cdot BDA^\top\|_{\mathrm{nuc}}$$
$$= \|A(\mathbf{I}_r - D^2)A^\top\|_{\mathrm{nuc}}$$
$$= r - \|D\|_\mathsf{F}^2 = r - \|U_t\|_\mathsf{F}^2$$
$$= r - \|(1-t)U + tV\|_\mathsf{F}^2$$
$$= r - (1-t)^2\|U\|_\mathsf{F}^2 - t^2\|V\|_\mathsf{F}^2 - 2t(1-t)\langle U, V\rangle$$
$$= r - (1-t)^2\|U\|_\mathsf{F}^2 - t^2\|V\|_\mathsf{F}^2 - t(1-t)\left(\|U\|_\mathsf{F}^2 + \|V\|_\mathsf{F}^2 - \|U - V\|_\mathsf{F}^2\right)$$
$$= r - (1-t)^2 r - t^2 r - t(1-t)\left(2r - \|U - V\|_\mathsf{F}^2\right)$$
$$= t(1-t)\|U - V\|_\mathsf{F}^2.$$

For (Term 2),

$$\|U_tU_t^\top - ((1-t)X + tY)\|_{\mathrm{nuc}} = \|((1-t)U + tV)((1-t)U + tV)^\top - (1-t)UU^\top - tVV^\top\|_{\mathrm{nuc}}$$
$$= \|-t(1-t)UU^\top - t(1-t)VV^\top + t(1-t)UV^\top + t(1-t)VU^\top\|_{\mathrm{nuc}}$$
$$= \|-t(1-t)(U-V)(U-V)^\top\|_{\mathrm{nuc}}$$
$$= t(1-t)\|U - V\|_\mathsf{F}^2.$$

Combining the two, then,
$$\min_{Z\in\mathcal{C}}\|Z - ((1-t)X + tY)\|_{\text{nuc}} \leq 2t(1-t)\|U - V\|_{\mathsf{F}}^2.$$

Next, note that the choice of $U$ and $V$ is not unique. Fixing any factorizations $X = UU^\top$ and $Y = VV^\top$, let $U^\top V = ADB^\top$ be a singular value decomposition, and let $\widetilde{V} = VBA^\top$. Then $Y = \widetilde{V}\widetilde{V}^\top$, and following the same steps as above we can calculate
$$\min_{Z\in\mathcal{C}}\|Z - ((1-t)X + tY)\|_{\text{nuc}} \leq 2t(1-t)\|U - \widetilde{V}\|_{\mathsf{F}}^2.$$

Furthermore,
$$\|U - \widetilde{V}\|_{\mathsf{F}}^2 = \|U\|_{\mathsf{F}}^2 + \|\widetilde{V}\|_{\mathsf{F}}^2 - 2\operatorname{trace}(U^\top \widetilde{V}) = 2r - 2\operatorname{trace}(U^\top VBA^\top)$$
$$= 2r - 2\operatorname{trace}(ADB^\top BA^\top) = 2r - 2\operatorname{trace}(D).$$

And,
$$\|X - Y\|_{\mathsf{F}}^2 = \|X\|_{\mathsf{F}}^2 + \|Y\|_{\mathsf{F}}^2 - 2\operatorname{trace}(XY) = 2r - 2\operatorname{trace}(UU^\top \widetilde{V}\widetilde{V}^\top)$$
$$= 2r - 2\|U^\top \widetilde{V}\|_{\mathsf{F}}^2 = 2r - 2\|D\|_{\mathsf{F}}^2 \geq 2r - 2\operatorname{trace}(D),$$

since $\|D\|_{\mathsf{F}}^2 = \sum_i (D_{ii})^2 \leq \sum_i D_{ii}$, as $0 \leq D_{ii} \leq 1$ for all $i$ since $U, V$ are both orthonormal matrices. Therefore, this proves that $\|U - \widetilde{V}\|_{\mathsf{F}}^2 \leq \|X - Y\|_{\mathsf{F}}^2$, and so
$$\min_{Z\in\mathcal{C}}\|Z - ((1-t)X + tY)\|_{\text{nuc}} \leq 2t(1-t)\|X - Y\|_{\mathsf{F}}^2.$$

Based on the curvature condition characterization (1) of the local concavity coefficients, we have therefore computed $\gamma_X(\mathcal{C}) \leq 2$, as desired. □